\documentclass[11pt,reqno]{amsart}

\usepackage{fixltx2e}                    
\usepackage{mathpazo}  
\usepackage[utf8]{inputenc}            
\usepackage{amsmath}                     
\usepackage{amssymb, latexsym, stmaryrd, amsthm, dsfont, amsfonts, amsbsy,amsthm, amsmath, mathrsfs}            
\usepackage{mathtools}                   
\usepackage{bm}                          
\usepackage{enumerate}                   

\usepackage{verbatim}                    
\usepackage{url}
\usepackage{lscape}                      
\usepackage{soul}

\usepackage[margin=0.85in]{geometry}

\usepackage{tikz} 

\usepackage{microtype}
\usepackage[all, knot]{xy}
        \xyoption{arc}
\makeatletter                            
\def\MT@register@subst@font{\MT@exp@one@n\MT@in@clist\font@name\MT@font@list
 \ifMT@inlist@\else\xdef\MT@font@list{\MT@font@list\font@name,}\fi}
\makeatother

\usepackage[pdftex,bookmarks,bookmarksnumbered,linktocpage,   
         colorlinks,linkcolor=blue,citecolor=blue]{hyperref}

\newcommand{\bit}{\begin{itemize}}    
\newcommand{\eit}{\end{itemize}}
\newcommand{\ben}{\begin{enumerate}}
\newcommand{\een}{\end{enumerate}}

\let\eroman\een

\newcommand{\benroman}{\ben[\normalfont (i)]}  
\newcommand{\benbullet}{\ben[\textbullet]}     
\let\ebullet\een
\newcommand{\bde}{\begin{description}}
\newcommand{\ede}{\end{description}}

\newcommand{\?}{\ensuremath{\mkern0.4\thinmuskip}}   
\let\models=\vDash                          

\let\leq=\leqslant
\let\nleq=\nleqslant

\let\geq=\geqslant

\let\epsilon=\varepsilon
\let\Lambda\varLambda
\let\Gamma\varGamma
\let\Delta\varDelta
\let\Lambda\varLambda
\let\Omega\varOmega
\let\Theta\varTheta
\let\Xi\varXi
\let\Pi\varPi
\let\Sigma\varSigma

\let\class=\mathsf                              
\let\oper=\mathbb                               

\bmdefine{\A}{A}                                
\bmdefine{\2}{2}
\bmdefine{\B}{B}
\bmdefine{\D}{D}
\bmdefine{\M}{M}                                
\bmdefine{\LLL}{L}                              
\bmdefine{\Fm}{Fm}                              
\bmdefine{\zerou}{[0{,}1]}
\bmdefine{\T}{T}                                

\newcommand{\integers}{{\mathbb Z}}
\newcommand{\naturals}{{\mathbb N}}

\newcommand{\HHH}{\oper{H}}

\newcommand{\PPP}{\oper{P}}

\newcommand{\SSS}{\oper{S}}
\newcommand{\III}{\oper{I}}

\theoremstyle{theorem}
\newtheorem{Theorem}{Theorem}[section]
\newtheorem{Proposition}[Theorem]{Proposition}
\newtheorem{Lemma}[Theorem]{Lemma}
\newtheorem{Corollary}[Theorem]{Corollary}

\newtheorem{Fan Theorem}[Theorem]{Brower's Fan Theorem}
\newtheorem{Matching Lemma}[Theorem]{Matching Lemma}
\newtheorem{Esakia Duality Theorem}[Theorem]{Esakia Duality Theorem}
\newtheorem{Claim}[Theorem]{Claim}

\theoremstyle{definition}
\newtheorem{law}[Theorem]{Definition}

\theoremstyle{remark}

\newtheorem{Remark}[Theorem]{Remark}


\begin{document}
\title[Profiniteness and representability of spectra of Heyting algebras]{Profiniteness and representability of spectra of Heyting algebras}

\author{G. Bezhanishvili, N. Bezhanishvili, T. Moraschini and M. Stronkowski}

\address{Guram Bezhanishvili: Department of Mathematical Sciences, New Mexico State University, Las Cruces NM 88003, USA}\email{guram@nmsu.edu}

\address{Nick Bezhanishvili: Institute for Logic, Language and Computation, University of Amsterdam, Postbus 94242, 1090GE Amsterdam, The Netherlands}\email{N.Bezhanishvili@uva.nl}

\address{Tommaso Moraschini: Department of Philosophy, University of Barcelona, Carrer de Montalegre $6$, $08001$, Barcelona, Spain}\email{tommaso.moraschini@ub.edu}

\address{Micha\l{} Stronkowski: Institute of Computer Science, Academy of Sciences of Czech Republic, Pod Vod\'arenskou v\v{e}\v{z}\'{i} $271/2$, $182$ $07$ Prague $8$, Czech Republic \newline {\rm and } \newline Faculty of Mathematics and Information Sciences,
Warsaw University of Technology, ul. Koszyko\-wa 75, 00-662
Warsaw, Poland}
\email{m.stronkowski@mini.pw.edu.pl}

\maketitle

\begin{abstract}
We prove that there exist profinite Heyting algebras that are not isomorphic to the profinite completion of any Heyting algebra.\ This resolves an open problem from $2009$.\ More generally, we characterize those varieties of Heyting algebras in which profinite algebras are isomorphic to profinite completions. It turns out that there exists largest such. We give different characterizations of this variety and show that it is finitely axiomatizable and locally finite. From this it follows that it is decidable whether in a finitely axiomatizable variety of Heyting algebras all profinite members are profinite completions. In addition, we introduce and characterize representable varieties of Heyting algebras, thus drawing connection to the classical problem of representing posets as prime spectra.
\end{abstract}

\tableofcontents

\section{Introduction}

An algebra is {\em profinite} if it is isomorphic to the inverse limit of an inverse system of finite algebras. This concept has its origin in the study of profinite groups. One of the basic results states that a group $\boldsymbol{G}$ is profinite if and only if it is a topological group whose topology is compact, Hausdorff, and zero-dimensional (i.e., a Stone topology). This result generalizes to many settings, including semigroups, monoids, rings, modules, distributive lattices, etc.~(see for instance \cite[Sec.~VI.2]{Jo82} or more recent \cite{CDJP08,SZ17}).  

The {\em profinite completion} $\widehat{\A}$ of an algebra $\A$ is the inverse limit of the inverse system of its finite homomorphic images. It follows from the definition that $\widehat{\A}$ is a profinite algebra. However, not every profinite algebra $\A$ is isomorphic to $\widehat{\B}$ for some algebra $\B$. In lattice theory, the problem of determining whether a profinite distributive lattice $\A$ is isomorphic to $\widehat{\B}$ for some distributive lattice $\B$ is related to Gr\"atzer's celebrated problem of representable posets. Gr\"atzer himself gave two formulations of this problem, for distributive lattices as for bounded distributive lattices (see \cite[Problems~33 and 34, p.~156]{Grtz71}). In this paper we will mainly concentrate on the bounded case. A poset $X$ is {\em representable} if it is isomorphic to the {\em prime spectrum} (the poset of prime filters) of a bounded distributive lattice. Gr\"atzer's problem asks for an internal characterization of representable posets. The same problem for prime spectra of commutative rings was posed by Kaplansky \cite[pp.~5--7]{Kaplansky1974}. In fact, bounded distributive lattices and commutative rings with unit give rise to the same prime spectra up to isomorphism (see, e.g., \cite[Thm.~1.1]{Pr94}).  Even more is true. It is a well-known result of Hochster \cite{Hochster1969} that the prime spectra of commutative rings with unit are precisely the spectral spaces. Cornish \cite{Cornish75} showed that these spaces form a category isomorphic to that of Priestley spaces, i.e., prime spectra of bounded distributive lattices \cite{Pr70,Pr72}.

It follows from Priestley duality that a poset $X$ is representable if and only if there is a Stone topology on $X$ such that the resulting ordered space is a Priestley space. Furthermore, in view of the work of Joyal \cite{Joy71} and Speed \cite{Spe72}, a poset $X$ is representable precisely when it is profinite. Thus, Gr\"atzer's problem asks for a description of profinite posets, or equivalently of the posets underlying Priestley spaces.  

A solution of Gr\"atzer's problem for linearly ordered sets was given independently by Balbes \cite[Thm.~9]{Bal71} (by lattice-theoretic means) and Lewis \cite[Thm.~3.1]{Lewis73} (by ring-theoretic means). It follows from their characterizations that representable chains are necessarily complete, whence a simple example of a chain that is not representable is given by the set of natural numbers with the usual order. 

Gr\"atzer's problem is connected to the problem of determining when a profinite bounded distributive lattice $\A$ is isomorphic to the profinite completion of some bounded distributive lattice $\B$, as we proceed to explain. Utilizing Priestley duality, it was shown in \cite[Thm.~4.4]{BeBe08} that $\A$ is a profinite bounded distributive lattice  if and only if $\A$ is isomorphic to the lattice $\textup{Up}(X)$ of all upsets (upward closed sets) of a poset $X$. On the other hand, by \cite[Thm.~5.3]{BeMo09}, $\A$ is isomorphic to $\widehat{\B}$ for some bounded distributive lattice $\B$ if and only if $\A$ is isomorphic to $\textup{Up}(X)$ for some representable poset $X$. Consequently, if the poset $X$ is not representable and $\A=\textup{Up}(X)$, then $\A$ is not isomorphic to $\widehat{\B}$ for any bounded distributive lattice $\B$ \cite[Cor.~5.4]{BeMo09}. This allows us to construct easily profinite bounded distributive lattices $\A$ that are not isomorphic to profinite completions of bounded distributive lattices; for instance, take $\A = \textup{Up}(X)$ where $X$ is any chain that is not complete.

In this paper we are interested in Heyting algebras. These are special bounded distributive lattices in which the meet operation has an adjoint, usually referred to as implication. Heyting algebras have been studied extensively as they have applications to different branches of mathematics, including:
\benbullet
\item {\bf Logic} (algebraic models of intuitionistic logic are Heyting algebras),
\item {\bf Topology} (the lattice of open sets of every topological space is a Heyting algebra), 
\item {\bf Point-free Topology} (each locale is a Heyting algebra), 
\item {\bf Domain Theory} (each continuous distributive lattice is a Heyting algebra), 
\item {\bf Topos Theory} (the subobject classifier of every topos is a Heyting algebra).
\item {\bf Lattice Theory/Universal Algebra} (every algebraic distributive lattice and hence the congruence lattice of every algebra in a congruence-distributive variety is a Heyting algebra).
\ebullet

Gr\"atzer's problem of representability was reformulated by Esakia for Heyting algebras in \cite[Appendix A.5]{Esakia-book85}. We call a poset $X$ {\em Esakia representable} if $X$ is isomorphic to the prime spectrum of some Heyting algebra. A poset $X$ is {\em image-finite} if the upset ${\uparrow}x \coloneqq \{y \in X :  x\leq y\}$ is finite for each $x\in X$. For a poset $X$, let
\[
X_{\textup{fin}} = \{x \in X :  {\uparrow}x\text{ is finite}\}
\]
be the image-finite subposet of $X$. By \cite[Thm.~3.6]{BeBe08}, a Heyting algebra $\A$ is profinite if and only if $\A$ is isomorphic to $\textup{Up}(X)$ for some image-finite poset $X$. By \cite[Thm.~5.8]{BeMo09}, a Heyting algebra $\A$ is isomorphic to $\widehat{\B}$ for some Heyting algebra $\B$ if and only if $\A$ is isomorphic to $\textup{Up}(X_{\textup{fin}})$ for some Esakia representable poset $X$. This implies that if there is an image-finite poset $X$ that is not isomorphic to $Y_{\textup{fin}}$ for any Esakia representable poset $Y$, then $\textup{Up}(X)$ is an example of a profinite Heyting algebra that is not isomorphic to the profinite completion of any Heyting algebra \cite[Cor.~5.9]{BeMo09}. Finding such an image-finite poset is nontrivial, and was left as an open problem in \cite{BeMo09}.

We will resolve this problem by providing many such examples. In fact, we will provide a characterization of all varieties of Heyting algebras whose profinite members are profinite completions of some Heyting algebras. This we do by proving that there exists a largest such variety.  We denote this variety by $\class{DHA}$ and term its members {\em diamond Heyting algebras}. To explain the terminology, $\class{DHA}$ is a subvariety of the variety $\class{CHA}$ of {\em cascade Heyting algebras} of Esakia \cite[Appendix A.9]{Esakia-book85}. We recall that $\class{CHA}$ is generated by {\em Boolean cascades}, i.e., algebras whose prime spectrum is a linear sum of finitely many antichains. If each antichain has at most two elements and we do not allow two ``back-to-back" two-element antichains, then the poset can be viewed as a linear sum of ``diamonds," and we refer to it as a {\em diamond sequence}. Let $\class{DHA}$ be the subvariety of $\class{CHA}$ generated by the Boolean cascades whose prime spectra are diamond sequences. This provides motivation for denoting this variety by $\class{DHA}$ and calling its members diamond Heyting algebras. We call prime spectra of diamond Heyting algebras {\em diamond systems}. Thus, diamond sequences are special diamond systems.

Our main result states that the profinite members of a variety $\class{V}$ of Heyting algebras are isomorphic to profinite completions of some Heyting algebras (which can always be chosen from $\class{V}$) if and only if $\class{V}$ is a subvariety of $\class{DHA}$ (Theorem~\ref{Thm:main}). Consequently, $\class{DHA}$ is the largest variety of Heyting algebras whose profinite members are profinite completions.

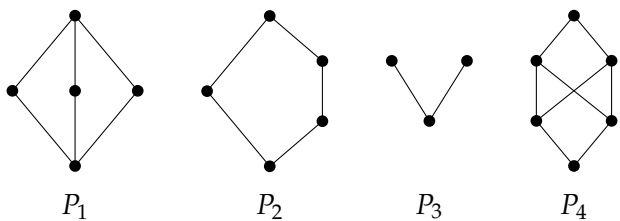
\begin{figure}[h]
\begin{tabular}{ccccccc}
\\
\begin{tikzpicture}
    \tikzstyle{point} = [shape=circle, thick, draw=black, fill=black , scale=0.35]
    \node[label=below:{$P_1$}]  (label) at (0,-0.1)  {};
    \node  (0) at (0,0) [point] {};
    \node (v1) at (-.833,1) [point] {};
    \node  (v2) at (0,1) [point] {};
    \node  (v3) at (.833,1) [point] {};
    \node (1) at (0,2) [point] {};

    \draw  (0) -- (v1) -- (1) -- (v2) -- (0) -- (v3) -- (1);
\end{tikzpicture}

&\,&
\begin{tikzpicture}
    \tikzstyle{point} = [shape=circle, thick, draw=black, fill=black , scale=0.35]
    \node[label=below:{$P_2$}]  (label) at (0,-0.1)  {};
    \node  (0) at (0,0) [point] {};
    \node (v1) at (-0.833,1) [point] {};
    \node  (v2) at (.7,0.6) [point] {};
    \node  (v3) at (.7,1.4) [point] {};
    \node (1) at (0,2) [point] {};

    \draw   (0) -- (v1) -- (1) --  (v3) -- (v2) -- (0);
\end{tikzpicture}
&\,&
\begin{tikzpicture}
    \tikzstyle{point} = [shape=circle, thick, draw=black, fill=black , scale=0.35]
    \node[label=below:{$P_3$}]  (label) at (0.5,-0.1)  {};
    \node at (0.5,0.0) {};  
    \node  (v1) at (0,1.4) [point] {};
    \node (v2) at (.5,.6) [point] {};
    \node  (v3) at (1,1.4) [point] {};

    \draw   (v1) -- (v2) -- (v3);
\end{tikzpicture}
&\,&
\begin{tikzpicture}
    \tikzstyle{point} = [shape=circle, thick, draw=black, fill=black , scale=0.35]
    \node[label=below:{$P_4$}]  (label) at (0.5,-0.1)  {};
    \node  (0) at (.5,0) [point] {};
    \node (v1) at (0,0.6) [point] {};
    \node  (v2) at (1,0.6) [point] {};
    \node  (v3) at (0,1.4) [point] {};
    \node  (v4) at (1,1.4) [point] {};
    \node (1) at (0.5,2) [point] {};

    \draw   (0) -- (v1) -- (v3) -- (1) -- (v4) -- (v2) -- (0);
    \draw   (v1) -- (v4);
    \draw   (v2) -- (v3);
\end{tikzpicture}
\end{tabular}
\caption{The posets $P_1$, $P_2$, $P_3$, and $P_4$.}
\label{fig:forbidden ordered-sets}
\end{figure}

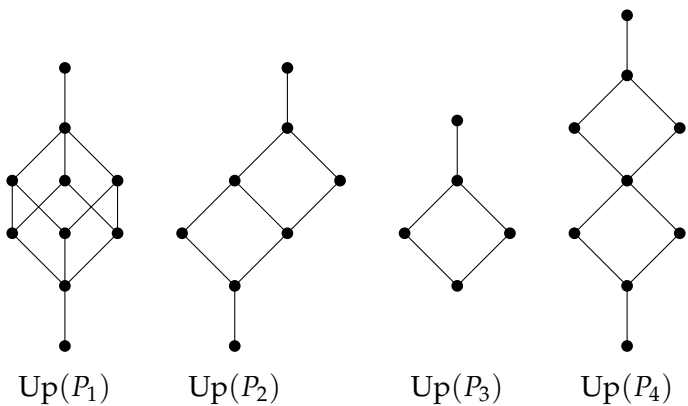
\begin{figure}[h]
\begin{tabular}{ccccccc}
\\
\begin{tikzpicture}
    \tikzstyle{point} = [shape=circle, thick, draw=black, fill=black , scale=0.35]
   
    \node[label=below:{$\textup{Up}(P_1)$}]  (label) at (0,-0.9)  {};
    \node  (-0) at (0,-.8) [point] {};
    \node  (0) at (0,0) [point] {};
    \node (v1') at (-.7,.7) [point] {};
    \node  (v2') at (0,.7) [point] {};
    \node  (v3') at (.7,.7) [point] {};
    \node (v1) at (-0.7,1.4) [point] {};
    \node  (v2) at (0,1.4) [point] {};
    \node  (v3) at (0.7,1.4) [point] {};
    \node (1) at (0,2.1) [point] {};
    \node  (+1) at (0,2.9) [point] {};
    \draw  (-0) -- (0) -- (v1') -- (v1) -- (1) -- (v3) -- (v3') -- (0) -- (v2') -- (v1);
    \draw   (v1') -- (v2) -- (v3') (v2') -- (v3) (v2) -- (1) -- (+1);
\end{tikzpicture}

&&
\begin{tikzpicture}
    \tikzstyle{point} = [shape=circle, thick, draw=black, fill=black , scale=0.35]
     \node[label=below:{$\textup{Up}(P_2)$}]  (label) at (0,-0.9)  {};
     
    \node  (-0) at (0,-.8) [point] {};
    \node  (0) at (0,0) [point] {};
    \node (v1) at (0.7,0.7) [point] {};
    \node  (v2) at (-.7,0.7) [point] {};
    \node  (v2') at (0,1.4) [point] {};
    \node  (v1') at (1.4,1.4) [point] {};
    \node (1) at (0.7,2.1) [point] {};
    \node (+1) at (0.7,2.9) [point] {};

    \draw   (-0)-- (0) -- (v1) -- (v1') --  (1) -- (v2') -- (v2)-- (0) (v1)--(v2') (1) -- (+1);
\end{tikzpicture}
&&
\begin{tikzpicture}
    \tikzstyle{point} = [shape=circle, thick, draw=black, fill=black , scale=0.35]
     \node[label=below:{$\textup{Up}(P_3)$}]  (label) at (0,-0.9)  {};
    \node at (0,-0.8) {};  
    \node at (0,-1.1){}; 
    \node  (0) at (0,0) [point] {};
    \node (v1) at (0.7,0.7) [point] {};
    \node  (v2) at (-.7,0.7) [point] {};
    \node (1) at (0,1.4) [point] {};
    \node (+1) at (0,2.2) [point] {};

    \draw   (0) -- (v1) -- (1) -- (v2) -- (0) (1) -- (+1);
\end{tikzpicture}
&&
\begin{tikzpicture}
    \tikzstyle{point} = [shape=circle, thick, draw=black, fill=black , scale=0.35]
     \node[label=below:{$\textup{Up}(P_4)$}]  (label) at (0,-0.9)  {};
    \node  (-0) at (0,-0.8) [point] {};
    \node  (0) at (0,0) [point] {};
    \node (v1) at (0.7,0.7) [point] {};
    \node  (v2) at (-.7,0.7) [point] {};
    \node (n) at (0,1.4) [point] {};
    \node (v1') at (0.7,2.1) [point] {};
    \node  (v2') at (-.7,2.1) [point] {};
    \node (1) at (0,2.8) [point] {};
    \node (+1) at (0,3.6) [point] {};

    \draw  (-0) -- (0) -- (v1) -- (n) -- (v1') -- (1) -- (v2') -- (n) -- (v2) -- (0) (1) -- (+1);
\end{tikzpicture}
\end{tabular}
\caption{The Heyting algebras of upsets of $P_1, P_2, P_3$, and $P_4$.}
\label{fig:forbidden Heyting algebras}
\end{figure}

Among various characterizations of diamond Heyting algebras that we give, one is of special interest. Let $P_{1}, P_{2}, P_{3}$, and $P_{4}$ be the posets shown in Figure \ref{fig:forbidden ordered-sets}.\footnote{In lattice theory it is customary to denote the posets $P_{1}$ and $P_{2}$ by $M_{3}$ and $N_{5}$.
\color{black}}
The Heyting algebras $\textup{Up}(P_{i})$ of upsets of these posets are shown in Figure \ref{fig:forbidden Heyting algebras}. We prove that the variety of diamond Heyting algebras is axiomatized by the Jankov formulas of $\textup{Up}(P_{i})$ for $i = 1, \dots, 4$.
\color{black}
Consequently, the problem of determining whether a variety of Heyting algebras is such that its profinite members are profinite completions of some Heyting algebras is decidable, both for varieties presented by a finite set of equations and for varieties presented by a finite number of finite algebras (Theorem~\ref{Thm:decidable}).

Our methods are based on a marriage of combinatorics of infinite posets, topology, and algebra. This is made possible by the use of Esakia duality \cite{Esakia74,Esakia-book85} that allows us to study Heyting algebras through the lenses of certain ordered topological spaces (known as Esakia spaces), which in turn are amenable to combinatorial principles such as Brouwer's Fan Theorem (see, e.g., \cite[Thm.~3.3.20]{Da04b}). Our main insight is a careful analysis of the structure of the  diamond systems that can be endowed with a Stone topology so that the resulting structures are Esakia spaces.

Some of our results have purely logical formulation. We recall that an \emph{intermediate logic} is a consistent axiomatic extension of  intuitionistic propositional logic (see, e.g., \cite[Sec.~2]{ChZa97}). These are exactly the logics that are situated between intuitionistic and classical logics and are algebraized by nontrivial varieties of Heyting algebras \cite{BP89}. Interemdiate logics algebraized by varieties of diamond Heyting algebra share many properties with the well-known G\"odel-Dummett logic (see, e.g., \cite[Sec.~5.1]{Dal86} and \cite[Sec.~4.2]{Haj98}). On the one hand, they are all finitely axiomatizable, locally tabular, and form a countable set (see Theorem~\ref{Thm:properties-of-DHAs}). On the other hand, they are all structurally complete \cite{Ry97} and have the infinite Beth definability property \cite{BlHoo06} (see Theorem~\ref{Thm:structural-ES}).

Finally, our results give a new insight on the representability problem as follows. Call a variety $\class{V}$ of Heyting algebras {\em representable} if $\textup{Up}(X)\in \class{V}$ implies that $X$ is Esakia representable. We give a full characterization of representable varieties by proving that a variety $\class{V}$ is representable if and only if $\class{V}\subseteq \class{DHA}$ and the depth of $\class{V}$ is bounded by some positive integer 
(see Theorem~\ref{Thm:representable-Esakia}). Consequently, every diamond system of bounded depth is Esakia representable. 
It follows that if a poset of finite depth validates the four axioms defining $\class{DHA}$, then it is Esakia representable. 
Since the axiomatization of $\class{DHA}$ is given by Jankov formulas of the Heyting algebras of the posets $P_1,\dots,P_4$ (see Figures~\ref{fig:forbidden ordered-sets} and~\ref{fig:forbidden Heyting algebras}), we obtain 
that if a poset of finite depth does not contain any of the posets $P_1,\dots,P_4$ as a forbidden configuration 
(a p-morphic image of an upset), then it is Esakia representable. In particular, every \textit{root system} (that is, a poset whose principal upsets are chains) \cite{NoGr15} of bounded depth is Esakia representable. It remains an interesting open problem to give a full characterization of Esakia representable posets of finite depth. A solution of this problem would shed further new light on the difficult problem of Esakia representability.

\section{Heyting algebras and Esakia spaces}

We start by recalling that a bounded distributive lattice $\A$ is a \emph{Heyting algebra} if $\land$ has an adjoint $\to$ given by
\[
a \land b \leq c \Longleftrightarrow a \leq b \to c
\]
for all $a, b, c \in A$.

It is well known (see, e.g., \cite{RaSi,BalDw,Esakia-book85}) that the class $\sf HA$ of Heyting algebras forms a variety (i.e., it is closed under homomorphic images, subalgebras, and direct products) and hence is equationally definable by Birkhoff's Theorem (see, e.g., \cite[Thm.~II.11.9]{BuSa00}).

\subsection{Esakia duality} The celebrated {\em Esakia duality} \cite{Esakia74,Esakia-book85} provides a representation of Heyting algebras by means of special ordered Stone spaces, and can be viewed as a restricted version of {\em Priestley duality} \cite{Pr70,Pr72} for bounded distributive lattices.

For a poset $X = \langle X,\leq \rangle$
\color{black}
and $U \subseteq X$, let ${\uparrow} U$ and ${\downarrow} U$ be the smallest upset and downset containing $U$, i.e.,
\begin{eqnarray*}
{\uparrow} U &=& \{x\in X :  \exists u\in U \mbox{ with } u\leq x\} \\
{\downarrow} U &=& \{x\in X :  \exists u\in U \mbox{ with } x\leq u\}
\end{eqnarray*}
We call $U$ an {\em upset} if $U={\uparrow}U$ and a {\em downset} if $U={\downarrow}U$. If $U = \{ x \}$, we simply write ${\uparrow} x$ and ${\downarrow} x$ instead of ${\uparrow} \{ x \}$ and ${\downarrow} \{ x\}$.

\begin{law}
An \emph{Esakia space} is a triple $X = \langle X, \tau, \leq \rangle$ where $\langle X, \tau \rangle$ is a Stone space (compact, Hausdorff, and zero-dimensional) and $\leq$ is a partial order on $X$ that is {\em continuous}; meaning that
\benroman
\item ${\uparrow} x$ is closed for all $x \in X$;
\item if $U\subseteq X$ is clopen, then ${\downarrow} U$ is clopen.
\eroman
\end{law}
\noindent Given an Esakia space $X$, we denote its underlying poset by $X$ as well.

\begin{Remark}
The partial order $\leq$ is continuous if and only if the corresponding map $\rho:X \to \mathcal V X$ from $X$ to the Vietoris space $\mathcal V X$, given by $\rho(x)={\uparrow}x$, is a well-defined continuous map (see \cite{Esakia74,Abr88,Kupke}).
\end{Remark}

As usual, we denote by $\textup{cl}(Y)$ the closure of a subset $Y$ of a topological space. Also, we denote the set of maximal elements of a poset $X$ by $\max X$. The following fundamental properties of Esakia spaces will be used repeatedly in the paper.

\begin{Proposition}\label{Prop:Esakia-tricks}
The following conditions hold for an Esakia space $X$.
\benroman
\item\label{Esakia-trick1} \emph{Priestley separation axiom:} for all $x, y \in X$ such that $x \nleq y$ there is a clopen upset $U$ such that $x \in U$ and $y \notin U$.
\item\label{item:chains-Esakia} \emph{Dedekind completeness:} Every nonempty chain $C \subseteq X$ has an infimum and a supremum in $X$.
\item\label{Esakia-trick3} If $C$ is a closed subset of $X$, then for every $x\in C$ there is $y \in \max C$ such that $x \leq y$.
\item\label{Esakia-trick5} For every $Y \subseteq X$, ${\uparrow}\textup{cl}(Y) = \textup{cl}({\uparrow}Y)$. Consequently, the closure of an upset is an upset and principal upsets are closed.
\item\label{Esakia-trick4} The downset of a closed set is closed. Consequently, principal downsets are closed.
\eroman
\end{Proposition}

\begin{proof}[Proof sketch.]
Conditions (\ref{Esakia-trick1}), (\ref{item:chains-Esakia}), and (\ref{Esakia-trick3}) are respectively \cite[Thms.~3.2.22(1), 3.2.19, and 3.2.1]{Esakia-book85}, while Conditions (\ref{Esakia-trick5}) and (\ref{Esakia-trick4}) follow from 
\cite[Thm.~3.1.2]{Esakia-book85}.
\end{proof}

In view of the above result, Esakia spaces satisfy Priestley separation axiom. In fact, Esakia spaces are those {\em Priestley spaces} (compact ordered spaces satisfying the Priestley separation axiom) in which the donwset of each (cl)open is (cl)open \cite[Thm.~3.1.2]{Esakia-book85}.

We recall that a {\em p-morphism} (or {\em bounded morphism}) between two posets $X$ and $Y$ is a map $f \colon X\to Y$ such that ${\uparrow}f(x)=f({\uparrow}x)$ for each $x\in X$.

\begin{law}
An \textit{Esakia morphism} between Esakia spaces is a continuous p-morphism. Let $\class{ES}$ be the category of Esakia spaces and Esakia morphisms between them.
\end{law}



\begin{Theorem} [Esakia duality]
$\class{HA}$ is dually equivalent to $\class{ES}$.
\end{Theorem}
We briefly describe the contravariant functors $(-)_*:\class{HA}\to\class{ES}$ and $(-)^*:\class{ES}\to\class{HA}$ establishing Esakia duality.
For this we first recall that for a poset $X$, the set $\text{Up}(X)$ of all upsets of $X$ forms a Heyting algebra where join and meet are set-theoretic union and intersection, and $\to$ is defined by
\[
U \to V = X \smallsetminus {\downarrow} (U \smallsetminus V) = \{ x\in X : {\uparrow}x\cap U\subseteq V \}.
\]

For a Heyting algebra $\A$, let $X_{\A}$ be the poset of prime filters of $\A$ ordered by inclusion. Define $\gamma_{\A} \colon \A\to\text{Up}(X_{\A})$ by
\[
\gamma_{\A}(a)=\{ x\in X_{\A} : a\in x \}.
\]

Then $\A_{\ast} = \langle X_{\A}, \tau, \subseteq \rangle$ is an Esakia space, where $\tau$ is the topology on $X_{\A}$ given by the subbasis
\[
\{ \gamma_\A(a) : a \in A \} \cup \{ \gamma_\A(a)^{c} : a \in A \}.
\]
If $\alpha \colon \A \to \B$ is a homomorphism between Heyting algebras, define $\alpha_{\ast} \colon \B_{\ast} \to \A_{\ast}$ by setting $\alpha_{\ast}(U) = \alpha^{-1}(U)$ for all $U \subseteq B_{\ast}$. Then $\alpha_{\ast}$ is an Esakia morphism, and this defines the contravariant functor $(-)_* \colon \class{HA}\to\class{ES}$.

For an Esakia space $X$, let $X^{\ast}$ be the subalgebra of $\textup{Up}(X)$ consisting of clopen upsets of $X$. If $f \colon X \to Y$ is an Esakia morphism between Esakia spaces, define $f^{\ast} \colon Y^{\ast} \to X^{\ast}$ by $f^{\ast}(U) = f^{-1}(U)$ for all $U \in Y^{\ast}$. Then $f^\ast$ is a homomorphism between the Heyting algebras $Y^{\ast}$ and $X^{\ast}$, and this defines the contravariant functor $(-)^* \colon \class{ES}\to\class{HA}$.

The topology of a finite Esakia space is discrete (since it is Hausdorff). Therefore, the full subcategory of $\class{ES}$ consisting of finite Esakia spaces is isomorphic to the category of finite posets and p-morphisms between them. Consequently, Esakia duality restricts to the following:

\begin{Theorem}[Finite Esakia duality]
The category of finite Heyting algebras is dually equivalent to the category of finite posets and p-morphisms between them.
\label{thm: FED}
\end{Theorem}

In addition, Stone spaces can be identified with the Esakia spaces whose underlying order is the identity. Under this identification, Esakia duality restricts to the famous {\em Stone duality} \cite{Sto36,Sto37} between the categories of Boolean algebras and Stone spaces.

To characterize homomorphic images and subalgebras of a Heyting algebra, we recall the notions of Esakia subspaces and Esakia quotients.


\begin{law}
Let $X$ be an Esakia space.
\begin{enumerate}[\normalfont (1)]
\item An \emph{Esakia subspace} ({\em E-subspace} for short)
is a closed upset of $X$ equipped with the subspace topology and the restriction of the order.
\item An \emph{Esakia equivalence} (or {\em E-partition} for short)
is an equivalence relation $R$ on $X$ satisfying for all $x, y, z \in X$:
\benroman
\item if $xRy$ and $y \leq z$, then $x\leq w$ and $wRz$ for some $w \in X$;
\item if $x {\not\!\! R} y$, then there is an {\em $R$-saturated} clopen $U$ (a union of equivalence classes of $R$) such that $x \in U$ and $y \notin U$.
\eroman
\end{enumerate}
\end{law}

\begin{Remark}
E-partitions are sometimes called correct partitions or bisimulation equivalences \cite{EsaGri77,Bez-PhD}.
\end{Remark}

It is well known that E-subspaces of Esakia spaces are Esakia spaces \cite[Lem.\ 3.4.11]{Esakia-book85}. In particular, closed subspaces of Stone spaces are Stone spaces.

If $R$ is an E-partition on $X$, we denote by $X / R$ the Esakia space consisting of the quotient space of $X$ with respect to $R$, equipped with the partial order $\leq_R$ defined as follows for every $x, y \in X$:
\begin{align*}
\llbracket x \rrbracket \leq_R \llbracket y\rrbracket \Longleftrightarrow& \text{ there are } x'\in \llbracket x\rrbracket \text{ and } y' \in \llbracket y\rrbracket \text{ such that } x' \leq y'.
\end{align*}
The map $x \mapsto \llbracket x \rrbracket$ is an Esakia morphism from $X$ to $X /R$, and for every Esakia morphism $f \colon X \rightarrow Y$, the kernel of $f$ is an E-partition of $X$ \cite[Thm.\ 2.3.9]{Bez-PhD}.

An algebra $\A$ is \textit{subdirectly irreducible} (SI for short) if the identity relation is completely meet-irreducible, and $\A$ is \emph{finitely subdirectly irreducible} (FSI for short) if the identity relation is meet-irreducible in the congruence lattice of $\A$ (see, e.g., \cite{Be11g,BuSa00}).

The following lemma collects some well-known consequences of Esakia duality.

\begin{Lemma}\label{Lem:correspondences}
Let $\A$ be a Heyting algebra.
\benroman
\item\label{Lem:correspondences:FSI} $\A$ is FSI iff its top element $1$ is \emph{prime} ($x \vee y = 1$ implies $x = 1$ or $y = 1$), which happens iff $\A_{\ast}$ is rooted (has a least element).
\item\label{Lem:correspondences:FactorAndSubspaces} The congruence lattice of $\A$ is dually isomorphic to the lattice of E-subspaces of $\A_{\ast}$.
\item\label{Lem:correspondences:SubalgebrasAndPartitions} The lattice of subalgebras of $\A$ is dually isomorphic to the lattice of E-partitions of $\A_{\ast}$.
\eroman
\end{Lemma}
\noindent When Condition (\ref{Lem:correspondences:FSI}) holds, $\{1\}$ is the least prime filter of $\A$  \cite[Prop.\ A.1.1]{Esakia-book85}. For Conditions (\ref{Lem:correspondences:FactorAndSubspaces}) and (\ref{Lem:correspondences:SubalgebrasAndPartitions}) see, e.g., \cite[Thm.\ 3.4.16]{Esakia-book85} and \cite[Thm.\ 2.3.9]{Bez-PhD}.

\subsection{Depth and width in Heyting algebras}

We next recall how to define finite depth and width for Esakia spaces and Heyting algebras, and prove that the set of elements of depth $\leq n$ (resp.~width $\leq n$) of an Esakia space forms an E-subspace.

\begin{law}
Let $X$ be a poset, $n$ a positive integer, and $x \in X$.
\begin{enumerate}[\normalfont (1)]
\item $X$ is said to have \emph{depth at most} $n$ if it does not contain any chain of $n+1$ elements.
\item $x \in X$ is said to be of {\em depth} $\leq n$ if the subposet ${\uparrow}x$ has depth $\leq n$.
\item $x \in X$ is said to be of {\em width} $\leq n$ if ${\uparrow} x$ does not contain any antichain of $n+1$ elements.
\item $X$ is said to have \emph{width at most} $n$ if every $x \in X$ is of width $\leq n$.
\end{enumerate}
\end{law}

\begin{law}
Let $\A$ be a Heyting algebra and $n$ a positive integer.
\begin{enumerate}[\normalfont (1)]
\item $\A$ has \emph{depth at most} $n$ if $\A_{\ast}$ has depth at most $n$. When this is the case, we write $d(\A)\leq n$.
\item $\A$ has \emph{depth exactly} $n$ if $d(\A)\leq n$ and $d(\A) \nleq n-1$.
\item $\A$ has \emph{width at most} $n$ if $\A_{\ast}$ has width at most $n$. When this is the case, we write $w(\A)\leq n$.
\item $\A$ has \emph{width exactly} $n$ if $w(\A)\leq n$ and $w(\A) \nleq n-1$.
\end{enumerate}
\end{law}

\noindent The notion of depth originates in \cite{Hosoi67}, while width was introduced for modal logics above $\class{K4}$ in \cite{Fine74}, 
and was adapted to Heyting algebras in \cite{Sob77}. 

Given a Heyting algebra $\A$ and a formula $\varphi$, we write $\A \vDash \varphi$ as a shorthand for $\A \vDash \varphi \thickapprox 1$. In this case, we say that $\varphi$ is \textit{valid} in $\A$. The following result is well known. 
\color{black}

\begin{Theorem}\label{Thm:d-w-varietal}
Let $\A$ be a Heyting algebra and $n$ a positive integer.
\benroman
\item\label{item:w:var:1} $d(\A)\leq n$  the formula $d_{n}$ is valid in $\A$,  where
\begin{align*}
d_{1} & = p_{1} \lor ( p_{1} \to 0)\\
d_{m+1} & = p_{m+1} \lor ( p_{m+1} \to d_{m}), \text{ for all }m  \geq 1.
\end{align*}
Consequently, the class $\class{D}_n$ of Heyting algebras of depth at most $n$ is a variety.
\item\label{item:w:var:2} $w(\A)\leq n$ iff the formula $w_{n}$ is valid in $\A$, where
\[
w_{n} = \bigvee_{i = 0}^{n}\big( p_{i} \to \bigvee_{j \ne i}p_{j}\big).
\]
Consequently, the class $\class{W}_n$ of Heyting algebras of width at most $n$ is a variety.
\eroman
\end{Theorem}

\begin{proof}[Proof sketch.]
Condition (\ref{item:w:var:1}) was essentially established in \cite{Ono71,Mak72} and Condition (\ref{item:w:var:2}) in \cite{Sob77}.
\end{proof}

Heyting algebras of width at most one are called \textit{G\"odel algebras} \cite[Sec.~4.2]{Haj98}. In view of the above result, they form a variety which, moreover, algebraizes the G\"odel-Dummett logic \cite{Dm59}.  Posets of width at most one have been called \textit{root systems} and coincide with disjoint unions of posets whose order duals are trees \cite{NoGr15}.


\begin{Proposition}\label{Prop:depth/width-E-subspaces}
Let $X$ be an Esakia space and $n$ a positive integer.
\benroman
\item\label{item-Prop:depth} The set of points of $X$ of depth $\leq n$ is an E-subspace of $X$.
\item\label{item-Prop:width}  The set of points of $X$ of width $\leq n$ is an E-subspace of $X$.
\eroman
\end{Proposition}

\begin{proof}
Condition (\ref{item-Prop:depth}) originates in \cite[Lem.~7]{Be00e}. Thus, we only give a proof of (\ref{item-Prop:width}). Let $X_{\leq n}$ be the set of points of $X$ of width $\leq n$. Since $x\leq y$ implies ${\uparrow}y\subseteq{\uparrow}x$, it is clear that $x\in X_{\leq n}$ implies $y\in X_{\leq n}$. Therefore, $X_{\leq n}$ is an upset of $X$. Thus, it is sufficient to show that $X_{\leq n}$ is closed in $X$. Since $X_{\leq n}$ is an upset, so is its closure $\textup{cl}(X_{\leq n})$ by Proposition \ref{Prop:Esakia-tricks}(\ref{Esakia-trick5}). Thus, $\textup{cl}(X_{\leq n})$ is an E-subspace of $X$. 

Let $\A$ be the Heyting algebra of clopen upsets of $\textup{cl}(X_{\leq n})$. For each $x\in X$, we have that ${\uparrow}x$ is an E-subspace of $X$ by Proposition~\ref{Prop:Esakia-tricks}(\ref{Esakia-trick5}). Let $\A_x$ be the Heyting algebra of clopen upsets of ${\uparrow}x$. Since for each $x\in X_{\leq n}$ we have ${\uparrow}x\subseteq X_{\leq n}$, we see that $\A_x\in \class{W}_n$.
Moreover, since $X_{\leq n}$ is dense in $\textup{cl}(X_{\leq n})$, the Heyting algebra $\A$ is a subalgebra of the product of $\{\A_x : x\in X_{\leq n}\}$ \cite[Appendix A.1]{Esakia-book85}. Therefore, $\A\in \class{W}_n$. Thus, the width of $\textup{cl}(X_{\leq n})$ is $\leq n$ by Theorem~\ref{Thm:d-w-varietal}(\ref{item:w:var:2}). Consequently,  $\textup{cl}(X_{\leq n}) = X_{\leq n}$, and hence $X_{\leq n}$ is closed in $X$.

Notice that, by interpreting $X_{\leq n}$ as the set of all elements of $X$ of depth $\leq n$, the above argument becomes a proof of (\ref{item-Prop:depth}).\end{proof}

\begin{Remark}
The above result fails for Priestley spaces. 
To see this, let $X$ be the one-point compactification of an infinite discrete space $Y$. We endow $X$ with a partial order by setting $x \leq y$ if and only if $x = y$ or $x$ is the unique limit point of $X$. It is well known and easy to check that the resulting ordered topological space is a Priestley space in which $\max X$ is dense. Thus, $\max X$ is the set of points of $X$ of both depth and width $1$. But its closure is $X$, whose depth is 2 and width is infinite. 
\end{Remark}

We next recall the following definition from the introduction.

\begin{law}
For an Esakia space $X$, let $X_{\textup{fin}}=\langle X_{\textup{fin}}, \leq \rangle$ be the subposet of $X$ where
\[
X_{\textup{fin}} = \{x\in X : {\uparrow}x \mbox{ is finite}\}.
\]
We call $X_{\textup{fin}}$ the {\em image-finite} part of $X$.
\end{law}

The next corollary to Proposition \ref{Prop:depth/width-E-subspaces} will be used throughout the paper.

\begin{Corollary}\label{Cor:finite-dw-topology}
Let $X$ be a poset of depth and width $\leq n$ for some positive integer $n$. If there exists an Esakia space whose image-finite part is $X$, then there is a Stone topology $\tau$ on $X$ such that $\langle X, \tau\rangle$ is an Esakia space.
\end{Corollary}

\begin{proof}
Suppose that $X$ is the image-finite part of an Esakia space $Y$. Let $Z$ be the set of elements of $Y$ of depth and width $\leq n$. By Proposition~\ref{Prop:depth/width-E-subspaces}, $Z$ is an E-subspace of $Y$. Therefore, it is sufficient to show that $X=Z$. By assumption, $X$ has depth and width $\leq n$. Since $X$ is an upset of $Y$, we see that $X \subseteq Z$. To prove the other inclusion, let $z \in Z$. Since $z$ has depth and width $\leq n$ in $Y$, it belongs to the image-finite part of $Y$, so $z \in X$. Thus, $X = Z$.
\end{proof}


\subsection{Profinite algebras and completions}

We recall that an \emph{inverse system} of Heyting algebras is a family $\{\A_i : i \in I \}$, indexed by a directed poset $I$, together with a family of homomorphisms $\alpha_{ij} \colon \A_j\to \A_i$ for $i\leq j$, satisfying:

\benroman
\item $\alpha_{kj}=\alpha_{ki}\circ \alpha_{ij}$ for all $k\leq i\leq j$;
\item $\alpha_{ii}$ is the identity homomorphism for each $i\in I$.
\eroman
The \emph{inverse limit} of such an inverse system is a Heyting algebra $\A$ together with a family $\alpha_i \colon \A\to \A_i$ of homomorphisms satisfying $\alpha_{ij}\circ \alpha_j=\alpha_i$ for each $i\leq j$, and having the following universal mapping property: For each Heyting algebra $\B$ and each family of homomorphisms $\beta_i \colon \B\to \A_i$ satisfying $\alpha_{ij}\circ \beta_j=\beta_i$ for each $i\leq j$, there is a unique homomorphism $\beta \colon \B\to \A$ such that $\alpha_i\circ \beta=\beta_i$ for each $i\in I$.

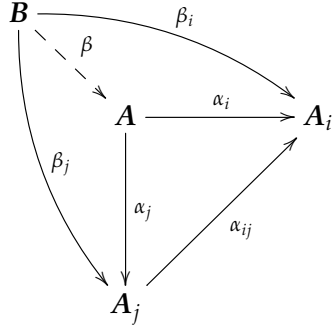
\[
\xymatrix{ \B  \ar@/_1.0pc/[dddr]^{\beta_j} \ar@/_-1.0pc/[rrrd]^{\beta_i} \ar@{-->}[dr]^\beta & & \\
& \A \ar[rr]^{\alpha_i} \ar[dd]^{\alpha_j} && \A_i \\
& & & \\
& \A_j \ar[uurr]_{\alpha_{ij}} &&}
\]

It is well known that $\A$ is isomorphic to the subalgebra of the product $\prod_{i\in I}\A_i$ consisting of
\[
\left\{ \vec{a} \in\prod_{i\in I}\A_i : \alpha_{ij}(\vec{a}(j))=\vec{a}(i) \mbox{ for each } j\geq i\right\}
\]
For inverse limits see, e.g., \cite[Ex.\ 13.13]{Adamek04}.

\begin{law}
We call a Heyting algebra $\A$ {\em profinite} if it is isomorphic to the inverse limit of an inverse system of finite Heyting algebras.
\end{law}

Let $\A$ be a Heyting algebra and $\theta$ a congruence on $\A$. We say $\theta$ has \emph{finite index} if the quotient algebra $\A/\theta$ is finite. Let $I$ be the set of congruences of $\A$ of finite index.
If $\theta\subseteq\phi$ are congruences of $\A$, then there is a canonical homomorphism $\pi_{\phi\theta} \colon \A/\theta\to \A/\phi$ given by
$\pi_{\phi\theta}(\llbracket a \rrbracket_\theta)= \llbracket a \rrbracket_\phi$. Furthermore, if $\theta$ and $\phi$ have finite index, then so does $\theta \cap \phi$ because $\A / (\theta \cap \phi)$ embeds into $\A / \theta \times \A / \phi$ via the map given by 
\[
\llbracket a \rrbracket_{\theta \cap \phi} \longmapsto \langle \llbracket a \rrbracket_{\theta}, \llbracket a \rrbracket_{\phi}\rangle,
\]
see, e.g., \cite[Ch.~II, Lem.~8.2]{BuSa00}. Therefore, $(I, \supseteq)$ is a directed poset, and $\{ \A/\theta : \theta \in I \}$ endowed with the homomorphisms $\pi_{\phi\theta}$ is an inverse system.


\begin{law}
The {\em profinite completion} $\widehat{\A}$ of a Heyting algebra $\A$ is the inverse limit of the inverse system described above.
\end{law}

Given a class $\class{K}$ of similar algebras, we denote by $\III(\class{K}), \HHH(\class{K}), \SSS(\class{K})$, and $\PPP(\class{K})$ the classes of isomorphic copies, homomorphic images, subalgebras, and direct products of elements of $\class{K}$, respectively. When $\class{K} = \{ \A \}$, we write $\III(\A)$ as a shorthand for $\III(\{ \A \})$. The same convention applies to $\HHH,\SSS$, and $\PPP$.

\begin{Remark}
If a Heyting algebra $\A$ is a profinite completion of some Heyting algebra $\B$, then $\A$ is also the profinite completion of some Heyting algebra $\boldsymbol{C}$ that belongs to the variety generated by $\A$; namely, $\HHH\SSS\PPP(\A)$.\ This is because $\A$ is also the profinite completion of $\boldsymbol{C} \coloneqq \B / \theta$, where $\theta$ is the intersection of all the congruences $\{ \theta_i : i \in I \}$ of $\B$ of finite index. Since $\B/ \theta_j \in \HHH(\A)$ for all $j \in I$ and $\B / \theta \in \III\SSS\PPP (\{ \B / \theta_i : \in I \})$, we obtain that $\boldsymbol{C} \in \III\SSS\PPP\HHH(\A) \subseteq \HHH\SSS\PPP(\A)$.
\end{Remark}
\color{black}

By \cite{BeHa16}, $\A$ is a profinite Heyting algebra if and only if there is a compact Hausdorff topology on $\A$ under which the Heyting operations are continuous. We will rely on the following characterization of profinite Heyting algebras and profinite completions of Heyting algebras.

\begin{Theorem}\label{Thm:folklore-profinite} Let $\A$ be a Heyting algebra and $X$ its Esakia dual.
\benroman
\item\label{folklore-profinite1} \cite[Thm.~3.6]{BeBe08} $\A$ is a profinite Heyting algebra if and only if there is an image-finite poset $Y$ such that $\A \cong \textup{Up}(Y)$.
\item\label{folklore-profinite2} \cite[Thm.~4.7]{BeCoMiMo06} $\widehat{\A}$ is isomorphic to $\textup{Up}(X_{\textup{fin}})$.
\eroman
\end{Theorem}

Clearly each profinite completion is a profinite Heyting algebra. However, as we will see, there are profinite Heyting algebras that are not isomorphic to the profinite completion of any Heyting algebra.
\color{black}


\section{Cascade Heyting algebras} 

In this section we recall cascade Heyting algebras and their dual description. We also give their axiomatization via Jankov formulas.

Let $X_1,\dots,X_n$ be Esakia spaces. We recall that the {\em linear sum} $X_1\oplus\cdots\oplus X_n$ is defined as the disjoint union of the spaces $X_1,\dots,X_n$, where the partial order is given by
\[
x\leq y \Longleftrightarrow (x,y\in X_i\text{ and }x\leq_i y\text{ for some }i\leq n)\text{ or }(x\in X_i\text{ and }y\in X_j\text{ for some }j < i \leq n).
\]
Figuratively speaking, we are forming a tower by putting $X_n,\dots,X_1$ on top of each other. We emphasize that $X_1$ is at the top of the tower and $X_n$ is at the bottom.
\color{black}
It is well known that a linear sum of Esakia spaces is again an Esakia space \cite[Prop.\ A.8.6]{Esakia-book85}. If $X$ is the linear sum of $X_1,\dots,X_n$ and $\A_1,\dots,\A_n$ are dual Heyting algebras to $X_1,\dots,X_n$, then the dual Heyting algebra $\A$ of $X$ is obtained by putting the algebras $\A_1,\dots,\A_n$ on top of each other and identifying the top element of $\A_i$ with the bottom element of $\A_{i+1}$ (see, e.g., \cite[Lem.~5.1]{MorWan19es}).

If each $X_i$ is a Stone space, then the linear sum $X=X_1\oplus\cdots\oplus X_n$ is a tower of antichains. Moreover, the dual of each $X_i$ is a Boolean algebra. Consequently, Esakia called a Heyting algebra whose dual Esakia space is a linear sum of Stone spaces a {\em Boolean cascade} \cite[Appendix A.9]{Esakia-book85} .

\begin{law}
Let $\class{CHA}$ be the variety of Heyting algebras generated by Boolean cascades. 
\end{law}

Following Esakia \cite[Appendix A.9]{Esakia-book85}, we call members of $\class{CHA}$ {\em cascade Heyting algebras}.

\begin{Theorem}\label{Thm:axiom:CHA-trick}
A Heyting algebra $\A$ is a cascade Heyting algebra iff it validates the \emph{weak Pierce law}:
\[
(p \to q) \lor (((q \to p) \to q) \to q).
\]
Thus, the logic ${\bf WPL}$ of the weak Pierce law is algebraized by the variety $\class{CHA}$ of cascade Heyting agebras.
\end{Theorem}

\begin{proof}[Proof sketch.]
This result was proved in \cite{Hosoi78,Esakia83}; see also \cite[Prop.~A.9.1(1)]{Esakia-book85}. A slightly different axiomatization of cascade Heyting algebras is given in \cite[Thm.~4]{Mardaev91}.
\end{proof}

An algebra is said to be \textit{locally finite} if its finitely generated subalgebras are finite. Accordingly, a variety is called \textit{locally finite} if all its members are locally finite. The following theorem originates with \cite{Mak81,Esakia83}; see also \cite[Prop.~A.9.1.(1)]{Esakia-book85}. A stronger result was proved by Kuznetsov \cite{Kuz74}; see, e.g., \cite[Thm.\ 2]{Citkin78a-paper}.

\begin{Theorem}\label{Thm:CHA-locally-finite}
$\class{CHA}$ is locally finite.
\end{Theorem}

In view of the above result and the fact that the class of Boolean cascades is closed under subalgebras, $\class{CHA}$ is generated by finite Boolean cascades. Hence, for many purposes, it is sufficient to work with linear sums of finite antichains.

Several useful characterizations of cascade Heyting algebras were given by Esakia \cite{Esakia-book85} and Mardaev \cite{Mardaev91}. The Mardaev characterization utilizes the technique of Jankov formulas \cite{Jankov63,Jankov68,Jankov69}. We recall that with each finite SI Heyting algebra $\A$, we can associate the \textit{Jankov formula}  $\mathcal{J}(\A)$ of $\A$, which axiomatizes the largest variety of Heyting algebras omitting the algebra $\A$. By Theorem~\ref{thm: FED} and Lemma~\ref{Lem:correspondences}(\ref{Lem:correspondences:FSI}), finite SI algebras correspond to finite rooted posets. Because of this, given a finite rooted poset $X$, we denote by $\mathcal{J}(X)$ the Jankov formula of the finite SI algebra $\textup{Up}(X)$.

\begin{Lemma}[Jankov's Lemma]\label{Lem:Jankov-trick-9}
Let $\A,\B$ be Heyting algebras with $\A$ finite and SI. Then $\B\models\mathcal{J}(\A)$ if and only if $\A\notin\mathbb{SH}(\B)$.
\label{lem: Jankov}
\end{Lemma}

In view of Lemma~\ref{Lem:correspondences}, Jankov's Lemma can be rephrased as follows.

\begin{Lemma}\label{lem:dual-jankov}
Let $X,Y$ be Esakia spaces with $X$ finite and rooted. Then $Y^{\ast}\models\mathcal{J}(X)$ if and only if $X$ is not a p-morphic image of an E-subspace of $Y$.
\end{Lemma}



\begin{law}[\protect{\cite[pg.\ 86]{Esakia-book85}}]
A poset $X$ is said to satisfy the \textit{three point rule} provided for every distinct $x, y, z \in X$, if $x$ and $y$ are incomparable, then $x\leq z$ implies $y \leq z$.
\end{law}

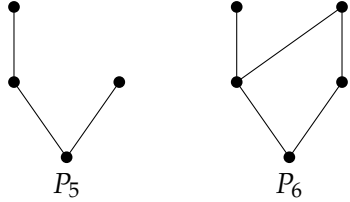
\begin{figure}[h]
\begin{tabular}{ccccccc}
\\
\begin{tikzpicture}
    \tikzstyle{point} = [shape=circle, thick, draw=black, fill=black , scale=0.35]
    \node[label=below:{$P_5$}]  (0) at (0,0) [point] {};
    \node (v1) at (-0.70,1) [point] {};
    \node  (v2) at (0.70,1) [point] {};
    \node (1) at (-0.70,2) [point] {};

    \draw  (0) -- (v1) -- (1)  (v2) -- (0) ;
\end{tikzpicture}

&&&&
\begin{tikzpicture}
    \tikzstyle{point} = [shape=circle, thick, draw=black, fill=black , scale=0.35]
    \node[label=below:{$P_6$}]  (0) at (0,0) [point] {};
    \node (v1) at (-0.7,1) [point] {};
    \node  (v2) at (0.7,1) [point] {};
    \node (1a) at (-0.7,2) [point] {};
    \node (1b) at (0.7,2) [point] {};

    \draw  (0) -- (v1) -- (1a) (v1) -- (1b) -- (v2) -- (0) ;
\end{tikzpicture}

\end{tabular}
\caption{The posets $P_5$ and $P_6$.}
\label{fig:P56}
\end{figure}

In the following theorem, $P_5$ and $P_6$ are the posets depicted in Figure \ref{fig:P56}.

\begin{Theorem}\label{Thm:Mardaev}
For a Heyting algebra $\A$, the following are equivalent.
\benroman
\item\label{item:Mardaev1} $\A$ is a cascade Heyting algebra;
\item\label{item:Mardaev2} ${\uparrow}x$ satisfies the three point rule for each $x\in\A_{\ast}$;
\item\label{item:Mardaev3} $\A$ validates the Jankov formulas $\mathcal{J}(P_2)$, $\mathcal{J}(P_5)$, and $\mathcal{J}(P_6)$.
\eroman
\end{Theorem}

\begin{proof}[Proof sketch.]
The equivalence of Conditions (\ref{item:Mardaev1}) and (\ref{item:Mardaev2}) goes back to \cite{Esakia83}; see also \cite[Prop.~A.9.2]{Esakia-book85}. 
The equivalence of Conditions (\ref{item:Mardaev1}) and (\ref{item:Mardaev3}) is a consequence of \cite[Thm.~5]{Mardaev91} and the fact that $\class{CHA}$ is generated by Boolean cascades.  
%
%
\end{proof}

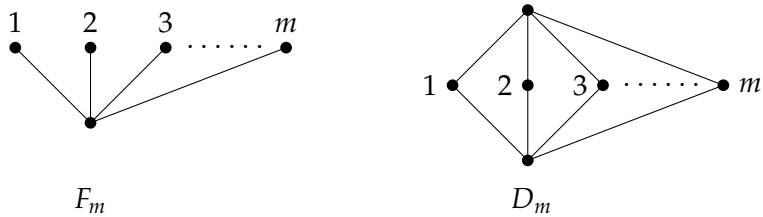
\begin{figure}[h]
\begin{tabular}{ccccccccc}
\\
\begin{tikzpicture}
     \tikzstyle{point} = [shape=circle, thick, draw=black, fill=black , scale=0.35]
    \node[label=below:{$F_m$}]  (0) at (-1.5,-0.6)  {};    
    \node[label=above:{$1$}]  (v1) at (-2.5,1) [point] {};
    \node (0) at (-1.5,0) [point] {};
    \node[label=above:{$3$}]  (v3) at (-0.5,1) [point] {};
    \node[label=above:{$2$}] (v2) at (-1.5,1) [point] {};
   \node (down-dots) at (0.32,1) {$\cdots\cdots$};
    \node[label=above:{$m$}]  (v4) at (1.1,1) [point] {};

    \draw   (v1) -- (0) -- (v3) (v4) -- (0) -- (v2);
\end{tikzpicture}

&&&&

\begin{tikzpicture}
     \tikzstyle{point} = [shape=circle, thick, draw=black, fill=black , scale=0.35]
    \node[label=below:{$D_m$}]  (0) at (-1.45,-0.1)  {}; 
    \node[label=left:{$1$}]  (v1) at (-2.5,1) [point] {};
    \node (0) at (-1.5,0) [point] {};
    \node[label=left:{$3$}]  (v3) at (-0.5,1) [point] {};
    \node[label=left:{$2$}] (v2) at (-1.5,1) [point] {};
   \node (down-dots) at (0.31,1) {$\cdots\cdots$};
    \node[label=right:{$m$}]  (v4) at (1.1,1) [point] {};
    \node (1) at (-1.5,2) [point] {};

    \draw (1) --  (v1) -- (0) -- (v3) -- (1) -- (v4) -- (0) -- (v2) -- (1);
\end{tikzpicture}

\end{tabular}
\caption{The posets $F_m$ and $D_m$.}
\label{fig:FDm}
\end{figure}

Let $\class{CHA}_n=\class{CHA}\cap \class{W}_n$. In other words, $\class{CHA}_n$ is the subvariety of $\class{CHA}$ consisting of cascade Heyting algebras of width $\leq n$. We conclude this section by showing that $\class{CHA}_n$ can be axiomatized by the Jankov formulas of the ``$(n+1)$-fork" $F_{n+1}$ and ``$(n+1)$-diamond" $D_{n+1}$. For every positive integer $m$, let $F_{m}$ and $D_m$ be the finite rooted posets depicted in Figure~\ref{fig:FDm}.


Our proof is based on the following two lemmas.

\begin{Lemma}\label{lem: 3.7}
If $\A\in\class{W}_n$, then $\A$ validates $\mathcal{J}(F_{n+1})$ and $\mathcal{J}(D_{n+1})$.
\end{Lemma}

\begin{proof}
Since $\A\in \class{W}_n$ and $\class{W}_n$ is a variety, each $\B\in\mathbb{SH}(\A)$ is also in $\class{W}_n$. Therefore, each such $\B$ has width $\leq n$. On the other hand, the width of both $\mathrm{Up}(F_{n+1})$ and $\mathrm{Up}(D_{n+1})$ is $n+1$. Therefore, neither belong to $\mathbb{SH}(\A)$. Thus, by Jankov's Lemma, $\A$ validates $\mathcal{J}(F_{n+1})$ and $\mathcal{J}(D_{n+1})$.
\end{proof}

\begin{Lemma}\label{lem: 3.8}
If $\A$ is a finite subdirectly irreducible cascade Heyting algebra and $\A\notin\class{W}_n$, then $\A$ refutes either $\mathcal{J}(F_{n+1})$ or $\mathcal{J}(D_{n+1})$.
\end{Lemma}

\begin{proof}
By finite Esakia duality, since $\A$ is finite and subdirectly irreducible, we may assume without loss of generality that $\A$ is $\mathrm{Up}(X)$ for some finite rooted poset $X$. In addition, since $\A$ is a finite cascade Heyting algebra and $X$ is rooted, $X$ satisfies the three point rule by Theorem~\ref{Thm:Mardaev}. It follows that $X$ is a linear sum of finite antichains. Since $\A\notin \class{W}_n$ and $X$ is rooted, there is an $(n+1)$-element antichain $C$ in $X$. Then $C$ is a subset of one of the antichains whose linear sum is $X$.  Let $x$ be an immediate predecessor of the elements of $C$ and let $Q={\uparrow}x$. Then $Q$ is an upset of $X$. 
 If $C \subseteq \max X$,  then $Q$ is isomorphic to $F_{n+1}$, so $\mathrm{Up}(F_{n+1})\in\mathbb{H}(\A)$. Therefore, by Jankov's Lemma, $\A\not\models \mathcal{J}(F_{n+1})$. Otherwise, define an equivalence relation $R$ on $Q$ by identifying all $y\notin C\cup\{x\}$. It is straightforward to check that $R$ is an E-partition of $Q$ and that $Q/R$ is isomorphic to $D_{n+1}$. Therefore, $\mathcal{J}(D_{n+1})\in\mathbb{SH}(\A)$. Thus, by Jankov's Lemma, $\A\not\models \mathcal{J}(D_{n+1})$.
\end{proof}

Putting together Lemmas~\ref{lem: 3.7} and~\ref{lem: 3.8} yields the following:

\begin{Theorem}\label{Lem:CHA-width2}
 $\class{CHA}_n$ is axiomatized over $\class{CHA}$ by the Jankov formulas $\mathcal{J}(F_{n+1})$ and $\mathcal{J}(D_{n+1})$.
\end{Theorem}

\begin{proof}
It follows from Lemmas~\ref{lem: 3.7} and~\ref{lem: 3.8} that $\class{CHA}_n$ and the subvariety of $\class{CHA}$ axiomatized by $\mathcal{J}(F_{n+1})$ and $\mathcal{J}(D_{n+1})$ have the same finite subdirectly irreducible members. The result follows since $\class{CHA}$ is locally finite by Theorem~\ref{Thm:CHA-locally-finite}. 
\end{proof}

\begin{figure}[h]
\begin{tabular}{ccccccccc}
\\
\begin{tikzpicture}
     \tikzstyle{point} = [shape=circle, thick, draw=black, fill=black , scale=0.35]
    \node at (0,0.0) {};  
    \node  (v1) at (-0.7,1) [point] {};
    \node (0) at (0,0) [point] {};
    \node  (v3) at (0.7,1) [point] {};
    \node (v2) at (0,2) [point] {};
    \node  (d3) at (1.4,2) [point] {};
    \node  (d1) at (-1.4,2) [point] {};

    \draw  (v2) -- (v1) -- (0) -- (v3)  -- (v2)  (d1) -- (v1) (v3) -- (d3);
\end{tikzpicture}

\end{tabular}
\caption{The poset $P_7$.}
\label{fig:W}
\end{figure}
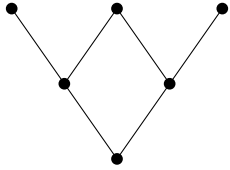

\begin{Remark}
The description of algebras of width $\leq n$ given in Theorem \ref{Lem:CHA-width2} cannot be extended to arbitrary Heyting algebras. For instance, the Heyting algebra of uspets of the poset $P_7$ depicted in Figure \ref{fig:W} validates $\mathcal{J}(F_{3})$ and $\mathcal{J}(D_{3})$ but is not of width $\leq 2$. Consequently, $\class{W}_n$ is not axiomatized over $\class{HA}$ by $\mathcal J(F_{n+1})$ and $\mathcal J(D_{n+1})$.
\end{Remark}

\section{Diamond systems and algebras}

 In this section we introduce the key concept of diamond systems and the corresponding variety $\class{DHA}$ of diamond Heyting algebras. We prove that $\class{DHA}$ is finitely axiomatizable by Jankov formulas. 


\begin{law}\label{def:diamond}
We call a poset $X$ a \emph{diamond system} if it satisfies the following conditions:
\benroman
\renewcommand{\labelenumi}{(\theenumi)}
\renewcommand{\theenumi}{D\arabic{enumi}}
\item\label{def:diamond1} ${\uparrow}x$ satisfies the three point rule for each $x\in X$;
\item\label{def:diamond2} $X$ has width at most two;
\item\label{def:diamond3} Principal upsets are upward directed in $X$;
\item\label{def:diamond4} For every $\bot, x, y, z, v, \top \in X$, if $\bot \leq x, y \leq z, v \leq \top$, there is $w \in X$ such that 
\[x, y \leq w \leq z, v.
\]
\eroman
If $X$ is in addition downward directed, then we call it a \textit{diamond sequence}.
\end{law}

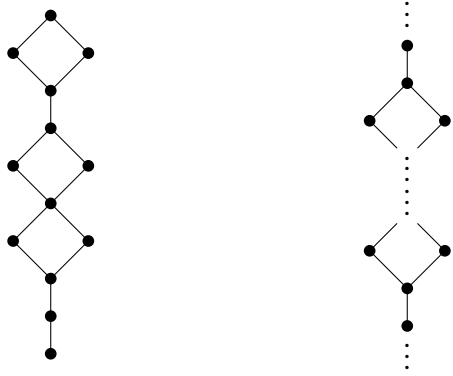
\begin{figure}[h]
\begin{tabular}{ccccccc}
\\
\begin{tikzpicture}
    \tikzstyle{point} = [shape=circle, thick, draw=black, fill=black , scale=0.35]
    \node  (b) at (0,0) [point] {};
    \node  (v1) at (0,0.5) [point] {};
    \node  (v2) at (0,1) [point] {};
    \node  (v3) at (-.5,1.5) [point] {};
    \node  (v3') at (.5,1.5) [point] {};
    \node  (v4) at (0,2) [point] {};
    \node  (v5) at (-.5,2.5) [point] {};
    \node  (v5') at (.5,2.5) [point] {};
    \node  (v6) at (0,3) [point] {};
    \node  (v7) at (0,3.5) [point] {};
    \node  (v8) at (-.5,4) [point] {};
    \node  (v8') at (.5,4) [point] {};
    \node  (t) at (0,4.5) [point] {};

    \draw   (b) -- (v1)  -- (v2) -- (v3) -- (v4) -- (v5) -- (v6) -- (v7) -- (v8) -- (t);
    \draw (v2)-- (v3') -- (v4) -- (v5') -- (v6)  (v7) -- (v8') -- (t);

    \node at (0,-0.2) {};   
\end{tikzpicture}

&\phantom{assadsadsadsads}&
\begin{tikzpicture}
    \tikzstyle{point} = [shape=circle, thick, draw=black, fill=black , scale=0.35]


   \node (v7) at (0,3.25) {$\vdots$};
   \node (v6) at (0,2.73) [point] {};
   \node (v5) at (0,2.23) [point] {};
   \node (v4) at (-.5,1.73) [point] {};
   \node (v4') at (.5,1.73)   [point] {};

   \node (v3) at (0,1.23)  {};
   \node (v2) at (0,1.2) {$\vdots$};
   \node (v2) at (0,0.75) {$\vdots$};
   \node (v1) at (0,.5)  {};
   \node (v0) at (0.5,0) [point] {};
   \node (v0') at (-0.5,0) [point] {};
   \node (v-1) at (0,-.5) [point] {};
   \node (v-2) at (0,-1) [point] {};
   \node (v-3) at (0,-1.3) {$\vdots$};

    \draw (v-2) -- (v-1) -- (v0) -- (v1) (v-1) -- (v0') -- (v1);
    \draw (v3) -- (v4) -- (v5) -- (v6) (v3) -- (v4') -- (v5);
\end{tikzpicture}
\end{tabular}
\caption{Two examples of diamond sequences.}
\label{Fig:diamond-sequences}
\end{figure}

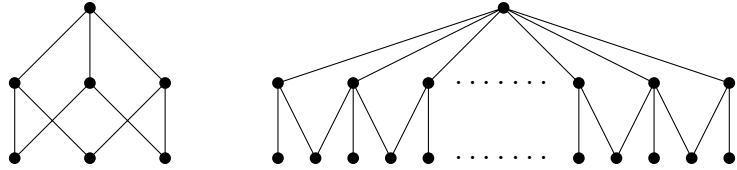
\begin{figure}[h]
\begin{tabular}{ccccccc}
\\
\begin{tikzpicture}
    \tikzstyle{point} = [shape=circle, thick, draw=black, fill=black , scale=0.35]
    \node  (t) at (-5.5,0) [point] {};
    \node  (x1) at (-6.5,-1) [point] {};
    \node  (x2) at (-5.5,-1) [point] {};
    \node  (x3) at (-4.5,-1) [point] {};
    \node  (y1) at (-6.5,-2) [point] {};
    \node  (y2) at (-5.5,-2) [point] {};
    \node  (y3) at (-4.5,-2) [point] {};

    \draw   (t) -- (x1)  -- (y1) -- (x2) -- (y3) -- (x3) -- (y2) -- (x1) (x2) -- (t) -- (x3);

   \node (t) at (0,0) [point] {};

   \node (x1)  at (-3,-1)  [point] {};
   \node (x2)  at (-2,-1)  [point] {};
   \node (x3)  at (-1,-1)  [point] {};
   \node (x4)  at (1,-1)  [point] {};
   \node (x5)  at (2,-1)  [point] {};
   \node (x6)  at (3,-1)  [point] {};

   \node (y1)  at (-3,-2)  [point] {};
   \node (y1')  at (-2.5,-2)  [point] {};
   \node (y2)  at (-2,-2)  [point] {};
   \node (y2')  at (-1.5,-2)  [point] {};
   \node (y3) at (-1,-2)  [point] {};
   \node (y4)  at (1,-2)  [point] {};
   \node (y4')  at (1.5,-2)  [point] {};
   \node (y5)  at (2,-2)  [point] {};
   \node (y5')  at (2.5,-2)  [point] {};
   \node (y6)  at (3,-2)  [point] {};

   \node (up-dots) at (0,-1) {$\cdots\cdot\cdots$};
   \node (down-dots) at (0,-2) {$\cdots\cdot\cdots$};

    \draw (t) -- (x1) -- (y1') -- (x2) -- (y2') -- (x3) -- (t) -- (x2);
    \draw (x1) -- (y1) (x2) -- (y2)  (x3) -- (y3);

    \draw (t) -- (x4) -- (y4') -- (x5) -- (y5') -- (x6) -- (t) -- (x5);
    \draw (x4) -- (y4) (x5) -- (y5)  (x6) -- (y6);

\end{tikzpicture}
\end{tabular}
\caption{An image-finite diamond system.}
\label{fig:: example of i-f almost forest}
\end{figure}

Notice that finite diamond sequences are simply linear sums of antichains of size $\leq 2$ in which we do not allow two ``back-to-back" antichains. Typical examples of diamond sequences are displayed in Figure \ref{Fig:diamond-sequences}, thus justifying the name. An image-finite diamond system that is not a diamond sequence is depicted in Figure \ref{fig:: example of i-f almost forest}. Observe that root systems coincide with diamond systems of width at most one.

Image-finite diamond systems admit a simple characterization in terms of linear sums.


\begin{Proposition}\label{Prop:shapes}
An image-finite poset $X$ is a diamond system if and only if for every $x \in X$ there are posets $Y_{1}, \dots, Y_{n}$ such that ${\uparrow}x$ is order-isomorphic to $Y_{1} \oplus \dots \oplus Y_{n}$, where $Y_{1}$ is a singleton, and for all $j > 1$, either  $Y_{j}$ is a singleton or isomorphic to $P_{3}$.
\end{Proposition}

\begin{proof}
The ``if'' part is straightforward. To prove the ``only if'' part, consider an image-finite diamond system $X$ and an element $x \in X$. Then ${\uparrow} x$ is finite since $X$ is image-finite. Therefore, by (\ref{def:diamond1}), ${\uparrow}x$ can be partitioned into a finite family of disjoint nonempty sets $X_{1}, \dots, X_{m}$ such that $X_{1} = \{ x \}$ and for every $y, z \geq x$,
\[
y \leq z \Longleftrightarrow \text{ either }y = z \text{ or }(y \in X_{i} \text{ and }z \in X_{j} \text{ for some }i < j).
\]
Notice that, as $X$ has width at most two by (\ref{def:diamond2}), each $X_{i}$ has at most two elements. Furthermore, if $X_{i}$ has two elements, then $i-1$ is positive, since $X_{1}$ is a singleton. Together with Condition (\ref{def:diamond4}), this implies that if $X_{i}$ has two elements, then $X_{i-1}$ exists and has one element. Thus, there are posets $Y_{1}, \dots, Y_{n}$ such that ${\uparrow}x$ is order-isomorphic to the sum $Y_{1} \oplus \dots \oplus Y_{n}$ where each $Y_{i}$ is either a singleton or isomorphic to $P_{3}$. To conclude the proof, it only remains to show that $Y_{1}$ is a singleton. But this is a consequence of the fact that $X$ is upward directed by (\ref{def:diamond3}).
\end{proof}

Notice that the assumption that $X$ is image-finite in the above result is essential since the chain of natural numbers (or the second poset in Figure \ref{Fig:diamond-sequences}) are diamond systems that do not satisfy Proposition \ref{Prop:shapes}.

\begin{law}
Let $\class{DHA}$ be the variety of Heyting algebras generated by the algebras of upsets of finite diamond sequences. The elements of $\class{DHA}$ will be called \textit{diamond Heyting algebras}.
\end{law}

Notice that, by definition, $\class{DHA}$ is a subvariety of $\class{CHA}$. Furthermore, since the variety $\class{GA}$ of G\"odel algebra is generated by the algebras of upsets of finite chains, which are diamond sequences, we obtain that $\class{GA}$ is a subvariety of $\class{DA}$.

Diamond systems and algebras are related by the following theorem, the hardest implication of which is (\ref{item:diamond4})$\Rightarrow$(\ref{item:diamond2}).

\begin{Theorem}\label{Thm:diamond}
The following conditions are equivalent for a Heyting algebra $\A$:
\benroman
\item\label{item:diamond1} $\A$ is a diamond Heyting algebra;
\item\label{item:diamond2} $\A_{\ast}$ is a diamond system;
\item\label{item:diamond3} Every principal upset in $\A_{\ast}$ is a diamond sequence;
\item\label{item:diamond4} $\A$ validates the Jankov formulas $\mathcal{J}(P_1), \mathcal{J}(P_2), \mathcal{J}(P_3)$, and $\mathcal{J}(P_4)$.
\eroman
\end{Theorem}

\begin{proof}
(\ref{item:diamond1})$\Rightarrow$(\ref{item:diamond4}): Notice that, by definition, diamond Heyting algebras are cascade Heyting algebras and have width $\leq 2$. Hence, by Theorem \ref{Thm:Mardaev} and Lemma \ref{Lem:CHA-width2}, they validate $\mathcal{J}(P_1)$ and $\mathcal{J}(P_2)$.  It only remains to prove that diamond Heyting algebras validate $\mathcal{J}(P_3)$ and $\mathcal{J}(P_4)$. In view of Lemma \ref{lem:dual-jankov}, it suffices to show that  $P_3$ and $P_4$  cannot be obtained as p-morphic images of upsets of finite diamond sequences. To this end, let $X$ be a finite diamond sequence. If $X$ is empty, we are done. Then suppose that $X \ne \emptyset$.  In this case, $X$ has a maximum, whence every nonempty p-morphic image of an upset of $X$ has also a maximum. It follows that $P_3$ cannot be a p-morphic image of an upset of $X$. The fact that $P_4$ is not a p-morphic image of an upset of $X$ follows from the fact that, by Proposition \ref{Prop:shapes}, $X$ can be viewed as linear sum $Y_{1} \oplus \dots \oplus Y_{n}$, where $Y_{1}$ is a singleton, and for all $j > 1$, either $Y_{j}$ is a singleton or isomorphic to $P_{3}$.

(\ref{item:diamond4})$\Rightarrow$(\ref{item:diamond2}): First we prove that $\A$ is a cascade Heyting algebra. To this end, observe that $P_3$ can be obtained as a p-morphic image of an upset of both $P_5$ and $P_6$. Consequently, if a variety of Heyting algebras contains either $\textup{Up}(P_5)$ or $\textup{Up}(P_6)$, it also contains $\textup{Up}(P_3)$. From the fact that $\A$ validates $\mathcal{J}(P_3)$ it follows that the variety generated by $\A$ omits $\textup{Up}(P_3)$, whence it also omits $\textup{Up}(P_5)$ and $\textup{Up}(P_6)$. We obtain that the Jankov formulas $\mathcal{J}(P_5)$ and $\mathcal{J}(P_6)$ are valid in $\A$. Recall also that, by assumption, $\A$ validates $\mathcal{J}(P_2)$. Hence, by Theorem \ref{Thm:Mardaev}, we conclude that $\A$ is a cascade Heyting algebra. From the same theorem it follows that principal upsets in $\A_\ast$ satisfy the three point rule. Thus, $\A_\ast$ satisfies Condition (\ref{def:diamond1}).

 We next  prove that $\A_\ast$ satisfies Condition (\ref{def:diamond2}), i.e.,  that it has width $\leq 2$. In view of Theorem~\ref{Lem:CHA-width2} it suffices to prove that $\A$ validates $\mathcal{J}(F_3)$ and $\mathcal{J}(D_3)$. Observe that $D_3 = P_1$ and that $\A$ validates $\mathcal{J}(P_1)$ by assumption. On the other hand, $P_3$ is a p-morphic image of $F_3$. It follows that if a variety of Heyting algebras contains $\textup{Up}(F_3)$, then it also contains $\textup{Up}(P_3)$. Since $\A$ validates $\mathcal{J}(P_3)$, the variety generated by $\A$ omits $\textup{Up}(P_3)$. It follows that the variety generated by $\A$  also omits  $\textup{Up}(F_3)$, whence $\A$ validates $\mathcal J(F_3)$. Thus, $\A_\ast$ has width $\leq 2$.

The next two claims conclude the proof that $\A_{\ast}$ is a diamond system.

\begin{Claim}\label{Claim:A-ast-upward}
Principal upsets are upward directed in $\A_{\ast}$.
\end{Claim}

\begin{proof}
Suppose not. Then there are $x, y, z \in A_{\ast}$ such that $x \leq y, z$ and ${\uparrow}y \cap {\uparrow}z = \emptyset$. As principal upsets are closed, we  can apply Proposition~\ref{Prop:Esakia-tricks}(\ref{Esakia-trick3}) to obtain that both ${\uparrow}y$ and ${\uparrow}z$ contain maximal elements. Since ${\uparrow}y \cap {\uparrow}z = \emptyset$, we may assume without loss of generality that these maximal elements are $y$ and $z$. Similarly, as principal downsets are closed, Proposition \ref{Prop:Esakia-tricks}(\ref{Esakia-trick3}) implies that there is an element that is maximal in ${\downarrow}z \cap {\downarrow}y$. We may assume without loss of generality that this element is $x$. Thus, $y$ and $z$ are maximal and $x$ is maximal among the lower bounds of $y$ and~$z$.

We show that
\begin{equation}\label{eq:diamond0}
{\uparrow}x = \{ x, y, z \}.
\end{equation}
The inclusion from right to left is clear. Suppose, with a view to contradiction, that the other inclusion fails. Then there is $a > x$ different from $y$ and $z$. Since $\A_{\ast}$ has width $\leq 2$ and $y$ and $z$ are incomparable (since ${\uparrow}y \cap {\uparrow}z = \emptyset$), $a$ is comparable either with $y$ or $z$. By symmetry we may assume that $a$ is comparable with $y$. Since $y$ is maximal and $a \ne y$, we get $a < y$. Because $x < a < y$ and $x$ is maximal among the lower bounds of $y$ and $z$, we have $a \nleq z$. As $z$ is maximal, this implies that $a$ and $z$ are incomparable. Since $x \leq a, y, z$ and $a \leq y$ and $a, y, z$ are different, the three point rule yields $a \leq z$, a contradiction.

Since principal upsets are closed in Esakia spaces,  ${\uparrow}x$ is the universe of an E-subspace $X$ of $\A_{\ast}$. From (\ref{eq:diamond0}) it follows that $X$ is isomorphic to $P_{3}$ endowed with the discrete topology. Hence, by Lemma \ref{Lem:correspondences}(\ref{Lem:correspondences:FactorAndSubspaces}), $\textup{Up}(P_3) \in \HHH(\A) \subseteq \class{DHA}$. But this contradicts the fact that $\A$ validates $\mathcal{J}(P_3)$.
\end{proof}

\begin{Claim}\label{Claim:A-ast-D4}
For every $\bot, x, y, z, v, \top \in A_{\ast}$, if $\bot \leq x, y \leq z, v \leq \top$, there is $w \in A_{\ast}$ such that $x, y \leq w \leq z, v$.
\end{Claim}

\begin{proof}
Suppose, with a view to contradiction, that there are $\bot, x, y, z, v, \top \in A_{\ast}$ such that $\bot \leq x, y \leq z, v \leq \top$ and ${\uparrow}x \cap {\uparrow}y \cap {\downarrow}z \cap {\downarrow}v = \emptyset$. Notice that $x$ and $y$ must be incomparable, otherwise $\{ x, y \} \cap {\uparrow}x \cap {\uparrow}y \cap {\downarrow}z \cap {\downarrow}v \ne \emptyset$, a contradiction. Similarly, $z$ is incomparable with $v$. Consequently, the elements $\bot, x, y, z, v, \top$ are all different.

Recall that $\bot$ is a lower bound of $x$ and $y$. Therefore, since principal donwsets are closed, Proposition \ref{Prop:Esakia-tricks}(\ref{Esakia-trick3}) implies the existence of a maximal lower bound of $x$ and $y$. We may assume without loss of generality that this maximal lower bound is $\bot$. Bearing this in mind, we will prove that
\begin{equation}\label{eq:diamond1}
{\uparrow}\bot = \{ \bot \} \cup {\uparrow}x \cup {\uparrow}y.
\end{equation}
The inclusion from right to left is clear. To prove the other inclusion,  consider an element $w > \bot$ different from $x$ and $y$. As $\A_{\ast}$ has width $\leq 2$, $w$ must be comparable either with $x$ or $y$. By symmetry we may assume that $w$ is comparable with $x$. If $w \geq x$, then $w \in \{ \bot \} \cup {\uparrow}x \cup {\uparrow}y$, as desired. Otherwise, $w < x$. We show this leads to a contradiction. Since $\bot$ is a maximal lower bound of $x$ and $y$ and $\bot < w \leq x$, we get $w \nleq y$. Moreover, $y \nleq w$, since $x$ and $y$ are incomparable and $w \leq x$. Thus, $x, y, w$ are different, $y$ and $w$ are incomparable and $w \leq x$. By the three point rule, $y \leq x$, a contradiction. This establishes (\ref{eq:diamond1}).

Furthermore, observe that
\begin{equation}\label{eq:diamond2}
{\uparrow}x \cup {\uparrow}y= \{ x, y, z, v  \} \cup ({\uparrow}z \cap {\uparrow}v).
\end{equation}
The inclusion from right to left is clear. Suppose, with a view to contradiction, that the other inclusion fails. Then there is $w \in A_{\ast} \smallsetminus \{ x, y, z, v \}$ such that $w \in {\uparrow}x \cup {\uparrow}y$ and $w \notin {\uparrow}z \cap {\uparrow}v$. Since $w \in {\uparrow}x \cup {\uparrow}y$, by symmetry we may assume that $x \leq w$. Similarly, as $w \notin {\uparrow}z \cap {\uparrow}v$, by symmetry we may assume that $z \nleq w$. Since $\A_{\ast}$ has width $\leq 2$, from $x \leq w, z, v$ and the fact that $z$ and $v$ are incomparable it follows that $w$ is comparable either with $z$ or $v$. We have two cases: either $w \in {\downarrow}z \cup {\downarrow}v$ or $w \notin {\downarrow}z \cup {\downarrow}v$. First suppose that $w \in {\downarrow}z \cup {\downarrow}v$. Since $w \notin \{ x, z, v \}$, either $x < w < z$ or $x < w < v$. By symmetry, we may assume that $x < w < z$. Since $x$ and $y$ are incomparable, $\bot \leq x, y$, and $x < w$, we can apply the three point rule to obtain $y \leq w$. It follows that $w \in {\uparrow}x \cap {\uparrow}y \cap {\downarrow}z$. Together with the assumption that ${\uparrow}x \cap {\uparrow}y \cap {\downarrow}z \cap {\downarrow}v$ is empty, this yields $w \nleq v$. Furthermore, from $w \leq z$ and the fact that $z$ and $v$ are incomparable it follows that $v \nleq w$. Hence, $w$ and $v$ are also incomparable. Therefore, since $x \leq w, v$ and $w < z$, we can apply the three point rule to obtain $v < z$, a contradiction.

Next  we consider the case where $w \notin {\downarrow}z \cup {\downarrow}v$. Together with the fact that $w$ is comparable either with $z$  or $v$ and  $z \nleq w$, this implies $v < w$. Since $z, v, w$ are distinct upper bounds of $x$ such that $z$ and $v$ are incomparable and $v \leq w$, we can apply the three point rule  to obtain $z \leq v$, a contradiction. This establishes (\ref{eq:diamond2}).

From (\ref{eq:diamond1}) and (\ref{eq:diamond2}) it follows that
\begin{equation}\label{eq:diamond3}
{\uparrow}\bot = \{ \bot, x, y, z, v \} \cup ({\uparrow}z \cap {\uparrow}v).
\end{equation}
Since ${\uparrow}\bot$ is a closed upset of $\A_{\ast}$, it is the universe of an E-subspace $X$ of $\A_{\ast}$. Similarly, ${\uparrow}z \cap {\uparrow}v$ is closed in $X$. In view of (\ref{eq:diamond3}), this implies that $\{ \bot, x, y, z, v \}$ is open in $X$. As $X$ is Hausdorff, we conclude that the points $\bot, x, y, z, v$ are isolated in $X$. Furthermore, $\{ \bot, x, y, z, v \}$ is closed because it is finite and $X$ is Hausdorff. By (\ref{eq:diamond3}) this yields that ${\uparrow}z \cap {\uparrow}v$ is open. Thus,
\begin{equation}\label{eq:diamond4}
\{ \bot \}, \{ x \}, \{ y \}, \{ z \}, \{ v \}, {\uparrow}z \cap {\uparrow}v\text{ are open in }X.
\end{equation}

We denote the elements of the poset $P_{4}$ as follows:
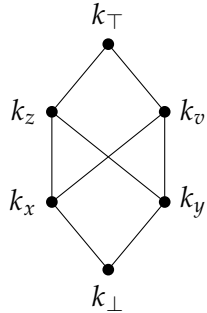
\[
\begin{tabular}{ccccccc}
\begin{tikzpicture}
    \tikzstyle{point} = [shape=circle, thick, draw=black, fill=black , scale=0.35]
    \node[label=below:{$k_{\bot}$}] (0) at (0,0) [point] {};
    \node[label=left:{$k_{x}$}] (v1) at (-0.75,0.9) [point] {};
    \node[label=right:{$k_{y}$}]  (v2) at (0.75,0.9) [point] {};
    \node[label=left:{$k_{z}$}]  (v3) at (-0.75,2.1) [point] {};
    \node[label=right:{$k_{v}$}]  (v4) at (0.75,2.1) [point] {};
    \node[label=above:{$k_{\top}$}] (1) at (0,3) [point] {};

    \draw   (0) -- (v1) -- (v3) -- (1) -- (v4) -- (v2) -- (0);
    \draw   (v1) -- (v4);
    \draw   (v2) -- (v3);
\end{tikzpicture}
\end{tabular}
\]
Bearing this in mind, consider the map $f \colon X \to P_{4}$ defined by
\[
f(a) = 
\begin{cases}
k_{a} & \text{ if } a \in \{ \bot, x, y, z, v \} \\
k_{\top} & \text{ if } a \in {\uparrow}z  \cap {\uparrow}v.
\end{cases}
\]
Notice that $f$ is a well-defined p-morphism by (\ref{eq:diamond3}). Moreover, when $P_{4}$ is endowed with the discrete topology, $f$ becomes an Esakia morphism. This is because $f^{-1}(k)$ is open in $X$ for every $k \in P_{4}$ by (\ref{eq:diamond4}). Thus, there is an E-partition $R$ of 
 $X$ such that $X / R \cong P_{4}$. Together with the fact that $X$ is an E-subspace of $\A_{\ast}$, this implies that 
 $\textup{Up}(P_{4}) \in \III\SSS\HHH(\A)$ by Lemma \ref{Lem:correspondences}(\ref{Lem:correspondences:FactorAndSubspaces}, \ref{Lem:correspondences:SubalgebrasAndPartitions}). But this contradicts the fact that $\A$ validates $\mathcal{J}(P_4)$.
\end{proof}

(\ref{item:diamond2})$\Rightarrow$(\ref{item:diamond3}): Immediate from the definitions.

(\ref{item:diamond3})$\Rightarrow$(\ref{item:diamond1}): For each $x \in A_\ast$, let $\A_x$ be the algebra dual to the E-subspace of $\A_\ast$ with the universe ${\uparrow}x$.\ From Esakia duality it follows that $\A$ is a subalgebra of the product of the family $\{ \A_x : x \in X \}$.\ Since varieties are closed under $\SSS$ and $\PPP$, it suffices to prove that each $\A_x$ is a diamond Heyting algebra. To this end, take some $x \in A_\ast$. Since $\A_\ast$ is a diamond system, the dual ${\uparrow}x$ of $\A_\ast$ is a diamond sequence. Accordingly, it satisfies the three point rule, whence $\A_x$ is a cascade Heyting algebra by Theorem \ref{Thm:Mardaev}. Hence, $\A_x$ is locally finite by Theorem \ref{Thm:CHA-locally-finite}. It follows that $\A_x$ belongs to the variety generated by its finite subalgebras. Therefore, to conclude the proof, it suffices to show that all finite subalgebras of $\A_x$ are diamond Heyting algebras. Consider one such algebra $\B$. The dual space $\B_\ast$ is a finite p-morphic image of the diamond sequence ${\uparrow}x$. Thus, $\B_\ast$ is a finite diamond sequence, whence $\B$ is a diamond Heyting algebra by definition.
\end{proof}

\begin{Corollary}\label{Cor:Axiomatization-jankov}
$\class{DHA}$ is axiomatized by $\mathcal{J}(P_1), \mathcal{J}(P_2), \mathcal{J}(P_3)$, and $\mathcal{J}(P_4)$.
\end{Corollary}

We close this section with the following corollary to Theorem \ref{Thm:diamond} which will be used in Section~\ref{Sec:6}.

\begin{Corollary}\label{Cor:profinite-come-from-diamond-systems}
Let $X$ be an image-finite poset. If $\textup{Up}(X)$ is a diamond algebra, then $X$ is a diamond system.
\end{Corollary}

\begin{proof}
As $\textup{Up}(X)$ is a diamond algebra, by Theorem \ref{Thm:diamond} the poset underlying $\textup{Up}(X)_{\ast}$ is a diamond system. It is well known that the map $\epsilon \colon X \to \textup{Up}(X)_{\ast}$, defined by 
\[
\epsilon(x) = \{ U \in \textup{Up}(X) : x \in U \} \text{ for all }x \in X,
\]
is a well-defined order embedding. Consequently, from the fact that $\textup{Up}(X)_{\ast}$ is a diamond system it follows that $X$ satisfies Conditions (\ref{def:diamond1}) and (\ref{def:diamond2}).

Suppose with a view to contradiction, that $X$ does not satisfy Condition (\ref{def:diamond3}). Then there are $x, y, z \in X$ such that $x \leq y, z$ and ${\uparrow}y \cap {\uparrow}z = \emptyset$. As $X$ is image-finite, ${\uparrow}x$ is finite. Thus, we may assume without loss of generality that $y$ and $z$ are maximal and that $x$ is maximal among the lower bounds of $y$ and $z$. Then we can repeat the proof of Claim \ref{Claim:A-ast-upward} and obtain the desired contradiction.

It only remains to prove that $X$ satisfies Condition (\ref{def:diamond4}). Suppose not. Then there are $\bot, x, y, z, v, \top \in X$ such that $\bot \leq x, y \leq z, v \leq \top$ and ${\uparrow}x \cap {\uparrow}y \cap {\downarrow}z \cap {\downarrow}v = \emptyset$. Since $X$ is image-finite, the upset ${\uparrow}\bot$ is finite, whence we may assume without loss of generality that $\bot$ is maximal among the lower bounds of $x$ and $y$. Thus, repeating the proof of Claim \ref{Claim:A-ast-D4} we obtain a contradiction, as desired.
\end{proof}

\section{Forbidden posets}

The aim of this section is to prove the following result.

\begin{Theorem}\label{Thm:first-half}
Let $\class{V}$ be a variety of Heyting algebras.\ If the profinite members of $\class{V}$ are profinite completions, then $\class{V} \subseteq \class{DHA}$.
\end{Theorem}

To this end, recall that a poset $\langle X, \leq \rangle$ is \emph{Esakia representable} if there exists a topology $\tau$ on $X$ such that $\langle X, \leq, \tau\rangle$ is an Esakia space. In order to establish Theorem \ref{Thm:first-half}, we require the following:

\begin{Proposition}\label{Prop:forbidden-posets}
If $P \in\{P_1,P_2,P_3,P_4\}$,  then there is a poset $X$ which is a p-morphic image of a disjoint union of copies of $P$ and is not Esakia representable.
\end{Proposition}

To derive Theorem \ref{Thm:first-half} from Proposition \ref{Prop:forbidden-posets} we recall the following well-known fact (the proof of which can essentially be found in \cite[Thm.\ 5.47(ii) \& 5.48]{BlRiVe01}).

\begin{Lemma}\label{Lem:preservation}
Let $\{ X_{i} \colon i \in I \} \cup \{ Y \}$ be a a family of posets and $\class{V}$ a variety of Heyting algebras such that $\{ \textup{Up}(X_{i}) \colon i \in I \} \subseteq \class{V}$. If $Y$ is a p-morphic image of the disjoint union of the various $X_{i}$, then $\textup{Up}(Y) \in \class{V}$.
\end{Lemma}

\begin{proof}[Proof of Theorem \ref{Thm:first-half}]
We reason by contraposition. Suppose $\class{V}$ is a variety of Heyting algebras  such that $\class{V} \nsubseteq \class{DHA}$. By Corollary \ref{Cor:Axiomatization-jankov}, $\class{DHA}$ is axiomatized by the Jankov formulas $\mathcal{J}(P_1), \mathcal{J}(P_2), \mathcal{J}(P_3)$, and $\mathcal{J}(P_4)$. Therefore, $\class{V} \nsubseteq \class{DHA}$ implies
\[
\class{V} \cap \{ \textup{Up}(P_1), \textup{Up}(P_2), \textup{Up}(P_3), \textup{Up}(P_4) \} \ne \emptyset.
\]
Consequently, there is $P \in\{P_1,P_2,P_3,P_4\}$ such that $\textup{Up}(P) \in \class{V}$. By Proposition \ref{Prop:forbidden-posets}, there is a p-morphic image $X$ of a disjoint union of copies of $P$ that is not Esakia representable. We have $\textup{Up}(X) \in \class{V}$ by Lemma \ref{Lem:preservation}. Furthermore, as $P$ has bounded depth and width, so does $X$. In particular, this guarantees that $X$ is image-finite, and hence $\textup{Up}(X)$ is a profinite member of $\class{V}$ by Theorem \ref{Thm:folklore-profinite}(\ref{folklore-profinite1}).

On the other hand, by Theorem \ref{Thm:folklore-profinite}(\ref{folklore-profinite2}) and Corollary \ref{Cor:finite-dw-topology}, the algebra $\textup{Up}(X)$ is a profinite completion if and only if $X$ is Esakia representable. As the latter is not the case, we conclude that $\textup{Up}(X)$ is not a profinite completion. Thus, $\class{V}$ has a profinite member that is not a profinite completion.
\end{proof}

The rest of this section is devoted to proving Proposition~\ref{Prop:forbidden-posets}.


%

\begin{law}
Let $X$ be a poset and $D,E\subseteq X$. We say that $(D,E)$ is a \emph{surjective matching} if the following conditions hold:
\benroman
\renewcommand{\labelenumi}{(\theenumi)}
\renewcommand{\theenumi}{M\arabic{enumi}}
\item \label{item:matching D cap E = emptyset} $D\cap E=\emptyset$;
\item\label{item:surjectivity on D} for every $x\in D$ there is $y\in E$ such that $y<x$;
\item\label{item:surjectivity on E} for every $y\in E$ there is $x\in D$ such that $y<x$;
\item\label{item:matching injectivity D} for every $x,z \in D$ and $y\in E$, if $y<x,z$ then $x=z$.
\eroman

A surjective matching $(D,E)$ is said to be a \textit{bijective matching} if it satisfies the following additional conditions:
\benroman
\renewcommand{\labelenumi}{(\theenumi)}
\renewcommand{\theenumi}{M\arabic{enumi}}
\setcounter{enumi}{4}
\item for every $x,y\in D\cup E$, if $y<x$ then $x\in D$ and $y\in E$;
\item\label{item:matching injectivity E} for every $x \in D$ and $y,z\in E$, if $y,z<x$ then $y=z$.
\eroman
In this case, $E \cup D$ becomes a disjoint union of two-element chains when endowed with the order inherited from $X$.
\end{law}

\begin{Lemma}[Matching Lemma]
\label{lem:: matching lemma}
Let $X$ be an Esakia space and $(D,E)$ a surjective matching in $X$. Let
\[
F = ({\uparrow}(D\cup E))\smallsetminus (D\cup E).
\]
If $E$ and $F$ are compact, then so is $D$.
\end{Lemma}

\begin{proof}
Suppose that $E$ is compact and consider an open cover $\mathcal K$ of $D$.  For $U\in \mathcal K$, define
\[
U' = {\downarrow}(U\smallsetminus (E\cup F)) \text{ and }{\mathcal K}' = \{U' : U\in{\mathcal K}\}.
\]
We will prove that ${\mathcal K}'$ is an open cover of $E$.

\begin{Claim} The following conditions hold:
\benroman
\item\label{item: U' is open} $U'$ is open;
\item\label{item: U' cap D = U cap D}  $U\cap D= U'\cap D$;
\item\label{item: x in U' iff y in U'} for every $x\in D$ and $y\in E$ with $y<x$ we have $x \in U'$ if and only if $y \in U'$.
\eroman
\end{Claim}

\begin{proof}
 (\ref{item: U' is open}): Since $X$ is Hausdorff and $E, F$ are compact, they are closed. Because $U$ is open, this implies that so is $U\smallsetminus (E\cup F)$. Lastly, since the downset of an open set is open in an Esakia space, we conclude that $U' = {\downarrow}(U\smallsetminus (E\cup F))$ is open, as desired.

(\ref{item: U' cap D = U cap D}):
By Condition (\ref{item:matching D cap E = emptyset}) and the definition of $F$, the intersection $(E\cup F)\cap D$ is empty. This yields
\[
U\cap D \subseteq  U\smallsetminus (E\cup F)\subseteq U'.
\]
Therefore, it remains to show that $U' \cap D \subseteq U$.

For this we first show that
\begin{equation}\label{eq:matching-1-xx}
({\uparrow}D)\smallsetminus D\subseteq F.
\end{equation}
Let $x \in ({\uparrow}D)\smallsetminus D$. We show that $x \notin E$. From $x \in ({\uparrow}D)\smallsetminus D$ it follows that there is $y \in D$ such that $y < x$.  If $x \in E$, by Condition (\ref{item:surjectivity on E}) there is $z \in D$ such that $x < z$. Since $y \in D$, by Condition (\ref{item:surjectivity on D}), $v < y$ for some $v \in E$. Now, $v \in E$ and $v < y, z \in D$. Therefore, $y = z$ by Condition (\ref{item:matching injectivity D}), contradicting $y < x < z$. This establishes $x \notin E$. Consequently,
\[
x \in ({\uparrow}D)\smallsetminus (D \cup E) \subseteq ({\uparrow}(D\cup E))\smallsetminus (D\cup E) = F,
\]
concluding the proof of (\ref{eq:matching-1-xx}).

Now suppose $x \in U' \cap D$. By (\ref{eq:matching-1-xx}),
\[
x \in U' =  {\downarrow}(U\smallsetminus (E\cup F)) \subseteq  {\downarrow}(U\smallsetminus  (({\uparrow}D)\smallsetminus D)).
\]
Therefore,  there is $y \in U \smallsetminus (({\uparrow}D)\smallsetminus D)$ such that $x \leq y$. From $x \in D$ and $x \leq y$, it follows  that $y \in {\uparrow}D$. Together with $y \in U \smallsetminus (({\uparrow}D)\smallsetminus D)$, this implies $y \in D$. Finally, in view of $x \in D$, by Condition (\ref{item:surjectivity on D}) there is $z \in E$ such that $z < x$. Consequently, $z \leq x, y \in D$. By Condition (\ref{item:matching injectivity D}) we obtain $x = y$, whence $x = y \in U$. Hence, $U' \cap D \subseteq U$.

(\ref{item: x in U' iff y in U'}): Suppose $x\in D$ and $y\in E$ with $y<x$.  Since $U'$ is a downset  and $y < x$, if $x\in U'$, then $y\in U'$. To prove the converse, suppose that $y\in U'$. Then there  is  $z\in U\smallsetminus (E\cup F)$ such that $y\leq z$. Notice that $z\in {\uparrow}E$ as  $z\geq y \in E$. Together with $z \notin F$, this implies $z\in D\cup E$. Furthermore, since $z \notin E$, we obtain  $z \in D$, whence $y \in E$, $x, z \in D$, and $y \leq x, z$. Thus, by Condition (\ref{item:matching injectivity D}),  $x = z \in U'$.
\end{proof}

Now, observe that $\mathcal{K}'$ is a family of open sets by Condition (\ref{item: U' is open}) of the Claim. We prove that it is a cover of $E$. Let $y\in E$. By Condition (\ref{item:surjectivity on E}), there is $x\in D$ such that $y<x$. Since $\mathcal K$ is a cover of $D$, by Condition (\ref{item: U' cap D = U cap D}) of the Claim, the same holds for  ${\mathcal K}'$. Thus, $x \in U'$ for some $U'\in \mathcal{K}'$. By Condition (\ref{item: x in U' iff y in U'}) of the Claim, we conclude that $y\in U'$. Hence, $\mathcal{K}'$ is an open cover of $E$.

Since $E$ is compact, $\mathcal{K}'$ has a finite subcover $\{U'_1,\ldots, U'_n\}$. We prove  that $\{U_1,\ldots, U_n\}$ is a cover of $D$. Let $x\in D$. By Condition (\ref{item:surjectivity on D}), there is  $y\in E$ such that $y<x$. Therefore,  $y\in U'_i$ for some $i$. Furthermore, $x\in U'_i$ by Condition (\ref{item: x in U' iff y in U'}) of the Claim.  Thus, applying Condition (\ref{item: U' cap D = U cap D}) of the Claim yields that $x\in U_i$. Consequently, $D$ is compact.
\end{proof}

We prove Proposition \ref{Prop:forbidden-posets} by contradiction.\ Once we have a candidate $\langle X, \leq \rangle$ for a poset which is not Esakia representable, we suppose that there is a topology $\tau$ on $X$ such that $\langle X, \leq, \tau \rangle$ is an Esakia space. Then we find a bijective matching $(D,E)$ in $X$ such that $E$ is compact and $D$ is not, a contradiction with the Matching Lemma.

For the compactness of $E$ we use the fact that in Esakia spaces principal downsets are closed, and hence compact. In the part of the proof regarding $P_1$ and $P_2$,
the set $E$ has the form ${\downarrow} x\cap {\downarrow }y$ for distinct elements $x$ and $y$. When dealing with $P_3$, we take $E$ of the form ${\downarrow}x\smallsetminus\{x\}$, where $x$ is an isolated point. The case of $P_4$ will be reduced to that of $P_{3}$.  The non-compact set $D$ will be constructed with the help of the following known result. We give a short proof for the sake of completeness.

\begin{Lemma}
\label{lem:non-compactification}
Let $X$ be a Hausdorff space and $Y$ an infinite subset of it. There is $y \in Y$ such that for every finite subset $Z$ of $Y$ the set $Y\smallsetminus(\{y\}\cup Z)$ is not compact in $X$.
\end{Lemma}

\begin{proof}
First suppose that $Y$ is not compact. In this case, subtracting from $Y$ any finite set gives a noncompact set. Thus, we can choose any $y \in Y$. Next  suppose that $Y$ is compact. Since $Y$ is infinite, this implies that $Y$ has a limit point $y$, i.e., an element $y\in Y$ such that every open neighbourhood $U$ of $y$ contains an element of $Y\smallsetminus\{y\}$. To conclude the proof, it suffices to show that $Y\smallsetminus\{y\}$ is not compact (as this implies that subtracting from $Y \smallsetminus \{ y \}$ any finite set produces a noncompact set). Since $X$ is Hausdorff, for every $x\in Y\smallsetminus \{y\}$, consider an open neighbourhood $U_{x}$ of $x$ and an open neighbourhood $V_{x}$ of $y$ such that $U_x\cap V_x=\emptyset$. Then $\{U_x : x\in Y\smallsetminus\{y\}\}$ is an open cover of $ Y\smallsetminus\{y\}$. Since $y$ is a limit point of $Y$, every intersection of finitely many sets of the form $V_x$ contains a point in $Y$. Hence, $\{U_x : x\in Y\smallsetminus\{y\}\}$ has no finite subcover of $Y\smallsetminus\{y\}$, whence $Y \smallsetminus \{ y \}$ is not compact.
\end{proof}

Lastly, we make use of the following observation.

\begin{Lemma}\label{Lem:Michal-new-lemma}
Let $X$ and $Y$ be posets and $\min X$ the set of minimal elements of $X$. If $X = {\uparrow}\!\?\min X$ and $Y \cong {\uparrow} x$ for every $x \in \min X$, then $X$ is a p-morphic image of a disjoint union of copies of $Y$. 
\end{Lemma}

\begin{proof}
For each $x \in \min X$, let $Y_{x}$ be a copy of $Y$ and let $f_{x} \colon Y_{x} \to {\uparrow}x$ be an order-isomorphism.
We let $Z$ be the disjoint union of the various $Y_{x}$. Then the union of the maps $f_{x}$ is a well-defined p-morphism $f \colon Z \to X$, which is onto since $X = {\uparrow}\!\? \min X$. Thus, $X$ is a p-morphic image of a disjoint union of copies of $Y$.
\end{proof}

\begin{proof}[Proof of Proposition \ref{Prop:forbidden-posets}] The proof is divided into four cases, corresponding to the posets $P_1, P_2, P_3$, and $P_4$.
\mbox{}

\subsubsection*{The case of $P_1$}

We construct a p-morphic image $X$ of a disjoint union of copies of $P_{1}$ that is not Esakia representable. For every three-element set  $\{k,m,n\}$ of natural numbers, consider a new element $\bot_{\{k,m,n\}}$. Also, consider another element $\top$  and  define
\[
X =\mathbb{N} \cup \{ \top \} \cup \{ \bot_{\{k,m,n\}} \colon k, m,n \text{ are distinct natural numbers} \};
\]
 for every $x, y \in X$, 
\begin{align*}
x \leq y \Longleftrightarrow &\text{ either } x = y \text{ or }y = \top\\ &\text{or }(x = \bot_{\{n,m,k\}} \text{ and }y \in \{ n, m, k \}\text{ for some }n, m, k \in \mathbb{N}).
\end{align*}

By Lemma \ref{Lem:Michal-new-lemma}, $X$ is a p-morphic image of a disjoint union of copies of $P_1$. Therefore, it only remains to prove that $X$ is not Esakia representable. Suppose the contrary, with a view to contradiction. Then there is an Esakia space $\langle X, \leq, \tau \rangle$. Choose a natural number $k$ and define $Y = \mathbb{N} \smallsetminus \{ k \}$. By Lemma \ref{lem:non-compactification}, there is $m \in \mathbb{N} \smallsetminus \{ k \}$ such that the set
\[
D \coloneqq \naturals\smallsetminus\{k,m\}
\]
is not compact. On the other hand, the set
\[
E \coloneqq \{\bot_{\{k,m,n\}} : n\in D\}={\downarrow}k\cap{\downarrow} m
\]
is the intersection of two principal downsets, so it is closed, and hence compact.
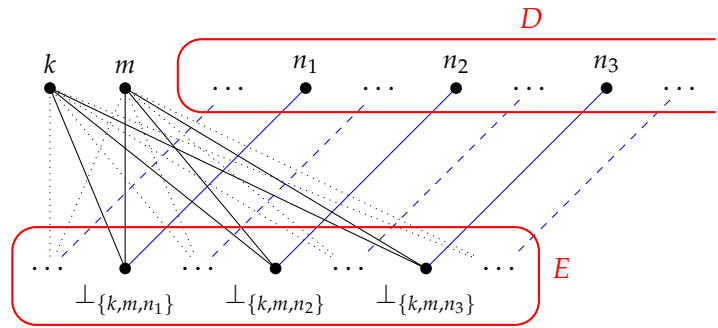
\begin{figure}[h]
\begin{tikzpicture}
    \tikzstyle{point} = [shape=circle, thick, draw=black, fill=black , scale=0.35]

    \node  (k) at (0,2.4) [point,label={\small$k$}] {};
    \node  (m) at (1,2.4) [point,label={\small$m$}] {};

    \node  (dots1) at (2.4,2.4) {$\dots$};
    \node  (n1) at (3.4,2.4) [point,label={\small$n_1$}] {};

    \node  (dots2) at (4.4,2.4) {$\dots$};
    \node  (n2) at (5.4,2.4) [point,label={\small$n_2$}] {};

    \node  (dots3) at (6.4,2.4)  {$\dots$};
    \node  (n3)   at (7.4,2.4) [point,label={\small$n_3$}] {};
    \node   (+inf) at (8.4,2.4)  {$\dots$};

    \node  (botdots1) at (0,0)  {$\dots$};
    \draw [dashed,color=blue] (botdots1) -- (dots1);
    \draw [dotted] (m) -- (botdots1) -- (k);

    \node  (bot1) at (1,0) [point,label=below:{\small$\bot_{\{k,m,n_1\}}$}] {};
    \draw [color=blue] (n1) -- (bot1);
    \draw (k) -- (bot1) -- (m);

    \node (botdots2) at (2,0)  {$\dots$};
    \draw[dashed,color=blue]  (botdots2) -- (dots2);
    \draw [dotted] (m) -- (botdots2) -- (k);

    \node (bot2) at (3,0) [point,label=below:{\small$\bot_{\{k,m,n_2\}}$}] {};
    \draw [color=blue] (n2) -- (bot2);
    \draw (k) -- (bot2) -- (m);

    \node (botdots3) at (4,0)  {$\dots$};
    \draw[dashed,color=blue] (dots3) -- (botdots3) ;
    \draw [dotted] (m) -- (botdots3) -- (k);

    \node  (bot3) at (5,0) [point,label=below:{\small$\bot_{\{k,m,n_3\}}$}] {};
    \draw [color=blue] (n3) -- (bot3);
    \draw (k) -- (bot3) -- (m);

    \node (botdots4) at (6,0)  {$\dots$};
    \draw[dashed,color=blue] (+inf) -- (botdots4);
    \draw [dotted] (m) -- (botdots4) -- (k);

    \draw[rounded corners=9pt, color=red, thick] (8.85,3.05) -- (1.7,3.05) -- (1.7,2.08) -- (8.85,2.08);

    \node [color=red] at (6.4,3.33) {$D$};

    \draw [rounded corners=9pt, thick, color=red] (-.5,-.75) rectangle (6.5,.55);

    \node [color=red] at (6.8,0) {$E$};

\end{tikzpicture}
\caption{The bijective matching $(D,E)$ in the poset $X$, constructed in the case of $P_1$.}
\label{fig:excluding M_3 b}
\end{figure}

\noindent Furthermore, the set
\[
F \coloneqq ({\uparrow}(D\cup E))\smallsetminus(D\cup E)=\{m,k,\top\}
\]
is finite. Consequently, as $X$ is a Hausdorff, $F$ is closed. In order to obtain a contradiction with the Matching Lemma, it only  remains to prove that $(D,E)$ is a surjective matching in $X$. But this follows immediately from the construction of $X$ and $(D,E)$, as illustrated in Figure \ref{fig:excluding M_3 b}.

\subsubsection*{The case of $P_2$}
We construct a p-morphic image $X$ of a disjoint union of copies of $P_{2}$ that is not Esakia representable. Let $\naturals' = \{1',2',3',\ldots\}$ be a disjoint copy of the set $\naturals$ of natural numbers. For every pair of distinct natural numbers $n$ and $k$, take a new element $\bot_{n,k}$. Consider two new elements $\top$ and $\sigma$ and define
\[
X = \naturals \cup \naturals' \cup \{ \top, \sigma \} \cup \{ \bot_{n,k} : n\text{ and } k \text{ are distinct natural numbers} \};
\]
for every $x, y \in X$,  
\begin{align*}
x \leq y \Longleftrightarrow &\text{ either }x = y \text{ or }y = \top\\
& \text{ or }(y = \sigma \text{ and }x \in X \smallsetminus (\naturals' \cup \{ \top \}))\\
&\text{ or }(x = \bot_{n,k} \text{ and }y \in \{ n , k' \}\text{ for some }n, k \in \mathbb{N}). 
\end{align*}
By Lemma \ref{Lem:Michal-new-lemma}, $X$ is a p-morphic image of a disjoint union of copies of $P_2$.  

Suppose, with a view to contradiction, that $X$ is Esakia representable. Then there is an Esakia space $\langle X, \leq, \tau\rangle$. By Lemma \ref{lem:non-compactification}, there exists a natural number $k$ such that $\naturals\smallsetminus\{k\}$ is not compact in $X$. Define
\[
D = \naturals\smallsetminus\{k\} \text{ and } E = \{\bot_{n,k} : n\in D\}={\downarrow}k'\cap{\downarrow}{\sigma}.
\]
Then $E$ is closed, hence compact. 
Moreover, as illustrated in Figure \ref{fig:excluding N5}, the set
\[
F \coloneqq ({\uparrow}(D\cup E))\smallsetminus(D\cup E)=\{\top,\sigma,k'\}
\]
is finite. Consequently, as $X$ is Hausdorff, $F$ is closed. In order to obtain a contradiction with Matching Lemma, it only remains to prove that $(D,E)$ is a surjective matching in $X$. But this is a direct consequence of the construction of $X$ and $(D,E)$, as shown in Figure \ref{fig:excluding N5}.

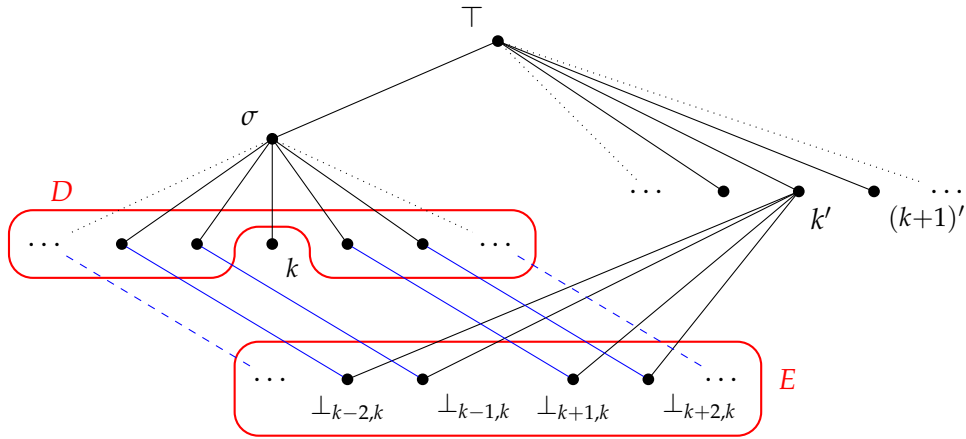
\begin{figure}[h]
\begin{tikzpicture}
    \tikzstyle{point} = [shape=circle, thick, draw=black, fill=black , scale=0.35]

   \node (top)   at (0,4.5) [point,label=above left:{$\top$}] {};
   \node (wedge) at (-3,3.2) [point,label=above left:{$\sigma$}] {};

    \node  (dots1) at (-6,1.8) {$\dots$};
    \node  (k-2) at (-5,1.8) [point] {};

    \node  (k-1) at (-4,1.8) [point] {};
    \node  (k) at (-3,1.8) [point,label=below right:{\small$k$}] {};

    \node  (k+1) at (-2,1.8)  [point] {};
    \node  (k+2)   at (-1,1.8) [point] {};
    \node  (dots2) at (0,1.8)  {$\dots$};

    \draw[thick, color=red, rounded corners=9pt] (0,1.35) -- (0.5,1.35) -- (0.5,2.25) -- (-6.5,2.25) -- (-6.5,1.35) -- (-3.5,1.35) -- (-3.5,2.03) -- (-2.5,2.03) -- (-2.5,1.35) -- (0,1.35);

    \node [color=red] at (-5.8,2.5) {$D$};

    \node  (bot dots1) at (-3,0) {$\dots$};
    \node  (bot k-2) at (-2,0) [point,label=below:{\small$\bot_{k{-}2,k}$}] {};

    \node  (bot k-1) at (-1,0) [point,label=below right:{\small$\bot_{k{-}1,k}$}] {};
    \node  (nothing) at (0,0) {};

    \node  (bot k+1) at (1,0)  [point,label=below:{\small$\bot_{k{+}1,k}$}] {};
    \node  (bot k+2)   at (2,0) [point,label=below right:{\small$\bot_{k{+}2,k}$}] {};
    \node  (bot dots2) at (3,0)  {$\dots$};

    \draw [rounded corners=9pt, thick, color=red] (-3.5,-.75) rectangle (3.5,.50);

    \node [color=red] at (3.85,0) {$E$};



     \node  (dots1') at (2,2.5) {$\dots$};

    \node  (k-1') at (3,2.5) [point] {};
    \node  (k') at (4,2.5) [point,label=290:{$k'$}] {};

    \node  (k+1') at (5,2.5)  [point,label=below right:{\small$(k{+}1)'$}] {};
    \node  (dots2') at (6,2.5)  {$\dots$};

    \draw[dashed] (bot dots1) (k') (bot dots2);
    \draw (bot k-2)--(k') -- (bot k+2);
    \draw (bot k-1)--(k') -- (bot k+1);



    \draw[dashed,color=blue] (dots1) -- (bot dots1);

    \draw[color=blue] (k-2) -- (bot k-2);
    \draw[color=blue] (k-1) -- (bot k-1);
    \draw[color=blue] (k+1) -- (bot k+1);
    \draw[color=blue] (k+2) -- (bot k+2);

    \draw[dashed,color=blue] (dots2) -- (bot dots2);


   \draw (wedge) -- (top) -- (k');
   \draw (k-1') -- (top) -- (k+1');
   \draw[dotted] (dots1') -- (top) -- (dots2');

   \draw (k-2) -- (wedge) -- (k+2);
   \draw (k-1) -- (wedge) -- (k+1);
   \draw (k) -- (wedge);
   \draw[dotted] (dots1) -- (wedge) -- (dots2);

\end{tikzpicture}
\caption{The bijective matching $(D,E)$ in the poset $X$, constructed in the case of $P_{2}$.}
\label{fig:excluding N5}
\end{figure}

\subsubsection*{The case of $P_3$}

We follow the same pattern as in the previous two cases. First, we construct a poset $X$ which is a p-morphic image of a disjoint union of copies of $P_{3}$. For every pair of integers $n$ and $k$ such that $n+1 < k$, consider a new element $\bot_{\{n,k\}}$ and define
\[
X = \integers \cup \{ \bot_{\{n,k\}} \colon n, k \in \integers \text{ and }n+1 < k \};
\]
for every $x, y \in X$,
\begin{align*}
x \leq y \Longleftrightarrow \text{ either }x = y
&\text{ or }(x = \bot_{\{n,k\}} \text{ and }y \in \{ n, k \} \text{ for some }n, k \in \mathbb{N}).
\end{align*}
By Lemma \ref{Lem:Michal-new-lemma}, $X$ is a p-morphic image of a disjoint union of copies of $P_3$. 

It only remains to show that $X$ is not Esakia representable. Suppose the contrary, with a view to contradiction. Then there exists an Esakia space $\langle X, \leq, \tau \rangle$. As in in the previous cases, we need to find a surjective matching $(D,E)$ in $X$ such that $E$ is compact, but $D$ is not. To this end, we rely on the following:

\begin{Theorem}[\protect{\cite[Thm.\ 8.5.4]{Semadeni}}]\label{Thm::every denumerable Stone space has an isolated point}
Every countable compact Hausdorff space has an isolated point.
\end{Theorem}

We shall construct a bijective matching in $X$ with the aid of the following:

\begin{Claim}
\label{claim:: finind D for F2}
There exists an integer $k$ that is isolated in $X$ and such that $\integers\smallsetminus \{k-1,k,k+1\}$ is not compact in $X$.
\end{Claim}

\begin{proof}
Notice that
$\integers$ is the set of maximal elements of $X$. Therefore, by Proposition \ref{Prop:depth/width-E-subspaces}, $\mathbb{Z}$ is closed in $X$. This means that the subspace $\integers$ of $X$ is a Stone space. Thus, we can apply Theorem \ref{Thm::every denumerable Stone space has an isolated point}, obtaining that $\integers$ has an isolated point $x$, i.e., a point $x$ such that $\{x\}$ is open in $\integers$. Since $\integers$ is compact and infinite, it must have a limit point. Hence, $\integers$ is the union of the two disjoint nonempty sets, respectively, of isolated points in $\integers$ and of limits points in $\integers$. As a consequence, there exist two consecutive integers such that one of them is a limit point in $\integers$ and the other is an isolated point in $\integers$ (otherwise either the set of isolated or  limit  points in $\integers$ would be empty).

Let $k$ be the isolated point in this pair. Then either $k-1$ or $k+1$ is a limit point in  $\integers$, whence the set $\integers\smallsetminus \{k{-}1,k,k{+}1\}$ is not compact in $\integers$,  and hence in $X$. To conclude the proof, it only remains to show that $k$ is isolated in $X$ too. Since $k$ is isolated in $\integers$, there exists an open set $U$ of $X$ such that $U\cap\integers=\{k\}$. Let
\[
V=X\smallsetminus {\downarrow}(X\smallsetminus U)=\{x \in X : {\uparrow}x \subseteq U\}.
\]
As downsets of closed sets are closed in Esakia spaces, $V$ is open in $X$. Furthermore, $V\smallsetminus\integers=\emptyset$  as  every element in $X\smallsetminus\integers$ has exactly two  successors in $\integers$ and one of them should be different from $k$. Together with the fact that $U\cap\integers=\{k\}$, this yields that $V=\{k\}$. Consequently,  $k$ is isolated in $X$.
\end{proof}

Let $k$ be as in Claim \ref{claim:: finind D for F2} and define
\[
D = \integers\smallsetminus\{k{-}1,k,k{+}1\} \text{ and } E = ({\downarrow}k) \smallsetminus\{k\}.
\]
Then $D$ is not compact and $E$ is closed, hence compact. From the definition of the order relation on $X$ it follows that
\[
F \coloneqq ({\uparrow}(D\cup E))\smallsetminus(D\cup E)=\{ k \}
\]
(see Figure \ref{fig:excluding F_2}). Thus, $F$ is finite,  and hence closed. As $(D,E)$ is a bijective matching in $X$ (see Figure \ref{fig:excluding F_2}),  this contradicts the Matching Lemma.

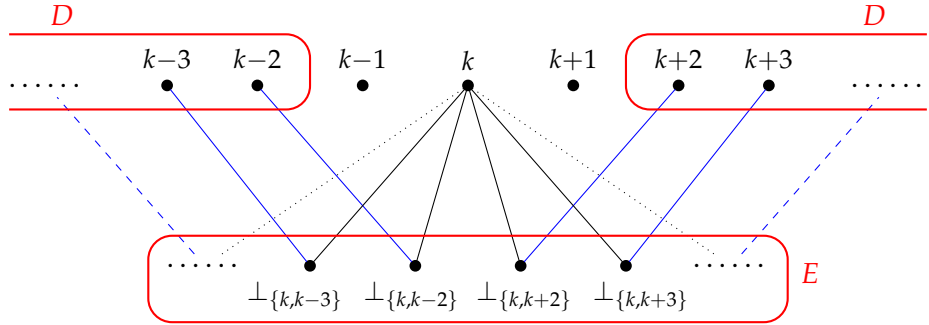
\begin{figure}[h]
\begin{tikzpicture}
    \tikzstyle{point} = [shape=circle, thick, draw=black, fill=black , scale=0.35]

    \node  (-inf) at (-5.6,2.4)  {$\dots\dots$};
    \node  (-3) at (-4,2.4) [point,label={\small $k{-}3$}] {};

    \node  (-2) at (-2.8,2.4) [point,label={\small $k{-}2$}] {};
    \node  (-1) at (-1.4,2.4) [point,label={\small $k{-}1$}] {};

    \node  (k) at (0,2.4) [point,label={\small $k$}] {};

    \node  (1) at (1.4,2.4) [point,label={\small $k{+}1$}] {};

    \node  (2) at (2.8,2.4)  [point,label={\small $k{+}2$}] {};
    \node  (3)   at (4,2.4) [point,label={\small $k{+}3$}] {};
    \node   (+inf) at (5.6,2.4)  {$\dots\dots$};

    \node  (-infbot) at (-3.5,0)  {$\dots\dots$};
    \node  (-3bot) at (-2.1,0) [point] {};

     \node at (-2.3,-0.4) {\small $\bot_{\{k,k{-}3\}}$}; 

    \node  (-2bot) at (-0.7,0) [point] {};

        \node at (-0.75,-0.4) {\small $\bot_{\{k,k{-}2\}}$}; 

    \node  (2bot) at (0.7,0)  [point,label] {};

    \node at (0.75,-0.4) {\small $\bot_{\{k,k{+}2\}}$}; 

    \node  (3bot)   at (2.1,0) [point] {};

        \node at (2.3,-0.4) {\small $\bot_{\{k,k{+}3\}}$}; 

    \node   (+infbot) at (3.5,0)  {$\dots\dots$};

    \draw [dashed,color=blue] (-infbot) -- (-inf);
    \draw [dashed,color=blue] (+infbot) -- (+inf);

    \draw [color=blue] (-3) -- (-3bot);
    \draw [color=blue] (-2) -- (-2bot);
    \draw [color=blue] (2) -- (2bot);
    \draw [color=blue] (3) -- (3bot);

   \draw [dotted] (-infbot) -- (k) -- (+infbot);
   \draw  (-3bot) -- (k) -- (3bot);
   \draw  (-2bot) -- (k) -- (2bot);

    \draw[rounded corners=9pt, color=red, thick] (6.1,3.08) -- (2.1,3.08) -- (2.1,2.08) -- (6.1,2.08);

    \draw[rounded corners=9pt, color=red, thick] (-6.1,3.08) -- (-2.1,3.08) -- (-2.1,2.08) -- (-6.1,2.08);

    \node [color=red] at (5.4,3.33) {$D$};
    \node [color=red] at (-5.4,3.33) {$D$};

    \draw [rounded corners=9pt, thick, color=red] (-4.25,-.75) rectangle (4.25,.4);

    \node [color=red] at (4.55,-.1) {$E$};

\end{tikzpicture}
\caption{The bijective matching $(D,E)$ in the poset $X$, constructed in the case of $P_{3}$.}
\label{fig:excluding F_2}
\end{figure}

\subsubsection*{The case of $P_4$}

Let $X$ be the poset defined in the case of $P_{3}$. We define a new poset $Y$ in which $X$ is a downset. Consider three new elements $\sigma, \tau$, and $\top$ and define
\[
Y = X\cup \{\sigma, \tau, \top\};
\]
for every $x, y \in Y$, 
\begin{align*}
x \leq y \Longleftrightarrow &\text{ either }x = y \text{ or } y=\top \text{ or } (x\in X \text{ and } y\in \{\sigma,\tau\})\\
&\text{ or }(x = \bot_{\{n,k\}} \text{ and }y \in \{ n, k \} \text{ for some }n, k \in \mathbb{N}).
\end{align*}
A portion of $X$ is depicted in Figure \ref{fig:big portion of the poset for excluding C2}.

\begin{figure}[h]

\begin{tikzpicture}
    \tikzstyle{point} = [shape=circle, thick, draw=black, fill=black , scale=0.35]

    \node  (top) at (0.5,4.8) [point] {};

    \node  (u) at (-0.4,4) [point] {};
    \node  (v) at (1.4,4) [point] {};

    \draw (u) -- (top) -- (v) ;

    \node   (-inf) at (-4,2.4)  {$\dots\dots$};
    \node  (-3) at (-3.2,2.4) [point] {}; 
    \node  (-2) at (-2.6,2.4) [point] {};
    \node (-1) at (-1.8,2.4) [point] {};
    \node  (0) at (-1,2.4) [point] {};
    \node  (k) at (-0.6,2.3)  {\small $k$};   
    \node  (1) at (2,2.4) [point] {};
    \node  (k+1) at (1.4,2.3)  {\small $k{+}1$};   
    \node  (2) at (2.8,2.4) [point] {};
    \node  (3) at (3.6,2.4) [point] {};
    \node  (4) at (4.4,2.4) [point] {};
    \node  (+inf) at (5.2,2.4)  {$\dots\dots$};

 \draw[dotted] (-inf) -- (u) -- (+inf) ;
    \draw (-3) -- (u) -- (4) ;
    \draw (-2) -- (u) -- (3) ;
    \draw (-1) -- (u) -- (2) ;
    \draw (0) -- (u) -- (1) ;

    \draw[dotted] (-inf) -- (v) -- (+inf) ;
    \draw (-3) -- (v) -- (4) ;
    \draw (-2) -- (v) -- (3) ;
    \draw (-1) -- (v) -- (2) ;
    \draw (0) -- (v) -- (1) ;

    \node (-inf!0) at (-2.3,0)  {$\dots$};
    \draw[dotted]   (-inf) -- (-inf!0) -- (0) ;
    \node (-3!0) at (-1.9,0) [point] {};
    \draw    (-3) -- (-3!0) -- (0) ;
    \node (-2!0) at (-1.6,0) [point] {};
    \draw    (-2) -- (-2!0) -- (0) ;
    \node (2!0) at (-1.3,0) [point] {};
    \draw    (2) -- (2!0) -- (0) ;
    \node (3!0) at (-1,0) [point] {};
    \draw    (3) -- (3!0) -- (0) ;
    \node (4!0) at (-0.7,0) [point] {};
    \draw    (4) -- (4!0) -- (0) ;
    \node  (+inf!0) at (-0.2,0)  {$\dots$};
    \draw[dotted]   (+inf) -- (+inf!0) -- (0) ;

    \node  (-inf!1) at (1.3,0)  {$\dots$};
    \draw[dotted]   (-inf) -- (-inf!1) -- (1) ;
    \node (-3!1) at (1.7,0) [point] {};
    \draw    (-3) -- (-3!1) -- (1) ;
    \node (-2!1) at (2,0) [point] {};
    \draw    (-2) -- (-2!1) -- (1) ;
    \node (-1!1) at (2.3,0) [point] {};
    \draw    (-1) -- (-1!1) -- (1) ;
    \node (3!1) at (2.6,0) [point] {};
    \draw    (3) -- (3!1) -- (1) ;
    \node (4!1) at (2.9,0) [point] {};
    \draw    (4) -- (4!1) -- (1) ;
    \node  (+inf!1) at (3.3,0)  {$\dots$};
    \draw[dotted]   (+inf) -- (+inf!1) -- (1) ;

\end{tikzpicture}
\caption{A fragment of the poset $Y$, constructed in the case of $P_{4}$.}
\label{fig:big portion of the poset for excluding C2}
\end{figure}
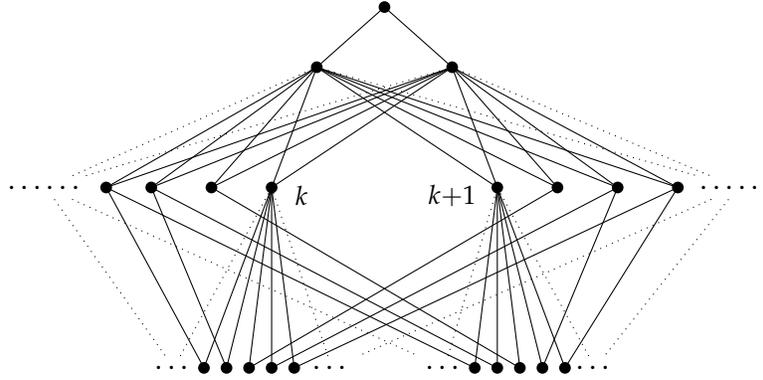

By Lemma \ref{Lem:Michal-new-lemma}, $X$ is a p-morphic image of a disjoint union of copies of $P_4$. Therefore, to conclude the proof, it suffices to show that $Y$ is not Esakia representable. As before, suppose with a view to contradiction, that there is an Esakia space $\langle Y, \leq, \tau \rangle$. We  prove that the ordered subspace $X$ of $Y$ is an Esakia space. For  this it suffices to show that $X$ is a clopen of $Y$ as, in this case, the conditions in the definition of an Esakia space follow \cite[Thm.\ 3.2.6]{Esakia-book85}. First notice that $X$ is open in $Y$ since $Y$ is Hausdorff and $Y\smallsetminus X$ is finite. Moreover, $X$ is closed in $Y$ since $X={\downarrow}{\sigma}\cap{\downarrow}{\tau}$. Consequently, $X$ is clopen in $Y$ and, hence, $X$ is an Esakia space. But this contradicts the fact that $X$ is not Esakia representable, which was established in the case of $P_{3}$.
\end{proof}

\section{Topologies for diamond systems}\label{Sec:6}

The aim of this section is to establish our main result:

\begin{Theorem}\label{Thm:main}
Let $\class{V}$ be a variety of Heyting algebras.\ The profinite members of $\class{V}$ are profinite completions if and only if $\class{V}$ is a subvariety of $\class{DHA}$.
\end{Theorem}

Notice that Theorem \ref{Thm:first-half} provides one half of the above result. To prove the other half, we require the following:

\begin{Proposition}
\label{prop:: topologies for almost forests}
Let $X$ be an image-finite diamond system. Then there exists an Esakia space $X^+$ such that its underlying poset is also a diamond system and  $X$ is the image-finite part of $X^+$.
\end{Proposition}

Proposition \ref{prop:: topologies for almost forests} implies Theorem \ref{Thm:main} as follows.

\begin{proof}[Proof of Theorem \ref{Thm:main}]
The ``only if'' part is Theorem \ref{Thm:first-half}. To prove the ``if'' part, suppose that $\class{V} \subseteq \class{DHA}$ and consider a profinite algebra $\A \in \class{V}$. By Theorem \ref{Thm:folklore-profinite}(\ref{folklore-profinite1}), there is an image-finite poset $X$ such that $\A \cong \textup{Up}(X)$. In particular, $\textup{Up}(X)$ is a diamond Heyting algebra, whence $X$ is a diamond system by Corollary \ref{Cor:profinite-come-from-diamond-systems}. In view of Proposition \ref{prop:: topologies for almost forests}, $X$ is the image-finite part of some Esakia space. Consequently, $\textup{Up}(X)$ is a profinite completion by Theorem \ref{Thm:folklore-profinite}(\ref{folklore-profinite2}). Since $\A \cong \textup{Up}(X)$, we conclude that $\A$ is a profinite completion.
\end{proof}

The rest of this section is devoted to proving Proposition \ref{prop:: topologies for almost forests}. 

\begin{proof}[Proof of Proposition \ref{prop:: topologies for almost forests}] Recall from Definition \ref{def:diamond} that diamond sequences are downward directed diamond systems. Let  ${\mathcal A}_\infty$ be the set of all infinite diamond sequences in $X$ that are maximal with respect to  inclusion.
With every $A\in\mathcal A_\infty$ we associate a new element $\bot_A$ (so that $A \ne B$ implies $\bot_A \ne \bot_B$). Let
\[
X^+ = X\cup \{\bot_A : A\in {\mathcal A}_\infty\}
\]
and extend the order on $X$ to that on $X^+$ by setting for $x,y\in X^+$, 
\begin{align*}
x\leq y \Longleftrightarrow  &\text{ either }x=y\text{ or }   x\leq y \text{ in }X \\
&\text{ or } x=\bot_A\text{ and } y\in A\text{ for some }A\in{\mathcal A}_\infty.
\end{align*}
The resulting poset $X^+$ is then obtained by adding a new lower bound $\bot_{A}$  for each maximal diamond sequence $A$. It follows that ${\uparrow}\bot_A = \{ \bot_A \} \cup A$. Notice that $X^+$ is a diamond system and $X$ is the image-finite part of $X^+$. In particular, if $X$ has finite depth, then $X = X^+$. 

To complete the proof, it suffices to show that $X^{+}$ can be endowed with a topology $\tau$ such that $\langle X^+, \tau, \leq \rangle$ is an Esakia space.  Notice that such a topology is not unique. For example, let $Y$ be an infinite poset of depth two with a maximum. We can turn $Y$ into an Esakia space by making any of its elements the unique limit point of the space. This results in different Esakia spaces with the same underlying poset $Y$.

To define $\tau$ on $X^+$, choose a maximal element $\top$ of $X$ (equiv.\ of $X^{+}$). Also, for every pair $x,y$ of incomparable elements in $X$ such that ${\downarrow}x\cap{\downarrow}y$ is nonempty, choose an immediate predecessor $\bot_{\{x,y\}}$ of $x$ and $y$ in $X$. The existence of $\bot_{\{x,y\}}$ follows from the fact that $X$ is image-finite and Proposition \ref{Prop:shapes}.

We define a topology $\tau$ on $X^+$ as follows. A set $U\subseteq X^+$ is open in $\tau$ if and only if it satisfies the following conditions: 

\vspace{1em}
\benroman
\renewcommand{\labelenumi}{(\theenumi)}
\renewcommand{\theenumi}{T\arabic{enumi}}
\item\label{item::top cond for top}
If $\top\in U$, then there is a finite set $Z\subseteq X$ such that $X^+\smallsetminus {\downarrow}Z\subseteq U$;

\item\label{item::top cond for down-set x} For every $x\in X$, if $x\in U$, then there is a finite set $Z\subseteq X\smallsetminus{\uparrow}x$ such that ${\downarrow}x\smallsetminus{\downarrow}Z\subseteq U$;

\item\label{item::top cond for bot x y}
For every pair of incomparable elements $x,y$ of $X$ such that ${\downarrow}x\cap{\downarrow}y\neq \emptyset$, if $\bot_{\{x,y\}}\in U$, then there is a finite set $Z\subseteq {\downarrow}x\cap{\downarrow}y$ such that $({\downarrow} x\cap{\downarrow}y)\smallsetminus{\downarrow}Z\subseteq U$;

\item\label{item::top cond for bot A}
For every $A\in\mathcal A_\infty$, if $\bot_A\in U$, then there is $x\in A$ such that ${\downarrow}x\subseteq U$.
\eroman

\begin{Remark}
 The following provides an intuitive meaning of $\tau$.  Postulate (\ref{item::top cond for top}) asserts that if $\top$ belongs to an open set $U$, then \textit{almost all} the space belongs to $U$. Similarly, (\ref{item::top cond for down-set x}) should be read as stating that if an element $x \in X$ belongs to $U$, then \textit{almost all} the downset ${\downarrow}x$ belongs to $U$. The remaining postulates can be explained along the same lines.
\end{Remark}

\begin{Claim}
The set $\tau$ is a topology on $X^+$.
\end{Claim}

\begin{proof}
That $\emptyset,X^+\in\tau$ and that $\tau$ is closed under arbitrary unions follow from the definition. Thus, it is enough to consider two sets $U_1, U_{2} \in \tau$ and verify that their intersection $U_1\cap U_2$ satisfies (\ref{item::top cond for top})--(\ref{item::top cond for bot A}).

First, suppose that $\top \in U_1\cap U_2$. Then there are finite $Z_1,Z_2\subseteq X$ such that $X^+\smallsetminus {\downarrow}Z_1\subseteq U_1$ and $X^+\smallsetminus {\downarrow}Z_2\subseteq U_2$. Therefore, $X^+\smallsetminus {\downarrow}(Z_1\cup Z_2)\subseteq U_1\cap U_2$. This shows that $U_{1} \cap U_{2}$ satisfies (\ref{item::top cond for top}). A similar argument shows that $U_{1} \cap U_{2}$ satisfies (\ref{item::top cond for down-set x}) and (\ref{item::top cond for bot x y}).

To show that $U_{1} \cap U_{2}$ satisfies (\ref{item::top cond for bot A}), suppose that $\bot_A\in U_1\cap U_2$ for some $A\in\mathcal A_\infty$. Then there are $x_1,x_2\in A$ such that ${\downarrow}x_1\subseteq U_1$ and ${\downarrow}x_2\subseteq U_2$. Since $A$ is a diamond sequence, there exists $x\in A$ such that $x\leq x_1,x_2$. Thus, ${\downarrow}x\subseteq U_1\cap U_2$, as desired.
\end{proof}

To prove that the ordered topological space  $X^+\coloneqq\langle X^+,\leq\,,\tau\rangle$ is an Esakia space, we show three facts: that downsets of open sets are open in $X^+$, that the space $X^+$ is compact, and that it satisfies the Priestley separation axiom. The verification of each of these is independent from the others.

\begin{Claim}\label{Fact:downsets-of-opens-are-opens}
Let $U$ be an open set in $X^+$. Then its downset ${\downarrow}U$ is also open in $X^+$.
\end{Claim}

\begin{proof}
We need to prove that ${\downarrow}U$ satisfies (\ref{item::top cond for top})--(\ref{item::top cond for bot A}). To prove (\ref{item::top cond for top}), suppose that $\top \in {\downarrow}U$. Since $\top$ is maximal, this yields $\top \in U$. As $U$ satisfies (\ref{item::top cond for top}), there is a finite set $Z$ such that $X \smallsetminus {\downarrow}Z \subseteq U$, whence also $X \smallsetminus {\downarrow}Z \subseteq {\downarrow}U$.

Condition (\ref{item::top cond for down-set x}) holds for every downset,  and in particular for ${\downarrow} U$. To prove (\ref{item::top cond for bot x y}), consider a pair of incomparable elements $x, y\in X$ such that ${\downarrow}x\cap{\downarrow}y\neq\emptyset$. Suppose that $\bot_{\{x,y\}}\in{\downarrow}U$. If $\bot_{\{x,y\}}\in U$, then the fact that ${\downarrow}U$ satisfies (\ref{item::top cond for bot x y}) follows from the assumption that $U$  satisfies (\ref{item::top cond for bot x y}) and the inclusion $U\subseteq {\downarrow}U$. On the other hand, if $\bot_{\{x,y\}}\notin U$, we have $\bot_{\{x,y\}}<u$ for some $u\in U$. Since $x$ and $y$ are immediate successors of $\bot_{\{ x, y \}}$, Proposition \ref{Prop:shapes} implies that $x\leq u$ or $y\leq u$. Thus, ${\downarrow}x\cap{\downarrow}y\subseteq{\downarrow}u\subseteq{\downarrow}U$.

It only remains to prove (\ref{item::top cond for bot A}). Suppose that $\bot_A\in{\downarrow}U$ for some $A\in \mathcal A_\infty$. If $\bot_A\in U$, then the fact that ${\downarrow}U$ satisfies (\ref{item::top cond for bot A}) follows from the assumption that $U$  satisfies (\ref{item::top cond for bot A}) and the inclusion $U\subseteq {\downarrow}U$. On the other hand, if $\bot_A \notin U$, there is $x\in U$ such that $\bot_A< x$. By the definition of the order  on $X^+$, we get $x\in A$. Therefore,  ${\downarrow}x\subseteq {\downarrow}U$.
\end{proof}


We next turn to proving  that $X^+$ is a compact space. Consider an open cover $\mathcal K$ of $X^+$. Our aim is to find a finite subcover of $\mathcal K$. The proof consists of a series of technical claims.

Let $U_\top$ be a member of $\mathcal K$ containing $\top$. By (\ref{item::top cond for top}), there exists a finite subset $Z$ of $X$ such that
\begin{equation}\label{Eq:the-definition-of-the-set-Z-}
X^+\smallsetminus U_{\top}\subseteq {\downarrow}  Z. 
\end{equation}
For every $z\in Z$, let $z^{\circ}$ be the unique maximal element in $X$ such that $z\leq z^{\circ}$. The existence and uniqueness of $z^{\circ}$ is guaranteed by the fact that $X$ is image-finite and Proposition \ref{Prop:shapes}. Define
\begin{equation}\label{Eq:definition-of-Y}
Y_1 = \{z^\circ : z\in Z\}, \,\, Y^+ = {\downarrow}Y_1, \text{ and } Y = Y^+\cap X,
\end{equation}
where ${\downarrow}Y_1$ is computed in $X^+$.

\begin{Claim}
\label{claim:: building finite branching S}
For every $x\in X$ there is a finite subset $\mathcal K_x$ of $\mathcal K$ and a finite subset $P_x$ of immediate predecessors of $x$  such that
\[
{\downarrow}x\smallsetminus {\downarrow}P_x\subseteq \bigcup {\mathcal K}_x.
\]
\end{Claim}

\begin{proof}
There is $U_x \in \mathcal K$ such that $x\in U_x$. By (\ref{item::top cond for down-set x}), there is a finite subset $Z_{x}$ of $X\smallsetminus {\uparrow}x$ such that ${\downarrow}x\smallsetminus {\downarrow}Z_{x}\subseteq U_x$. Clearly we may assume that ${\downarrow}x\cap{\downarrow}z\neq \emptyset$ for every $z\in Z_{x}$. Define
\[
Z_1^{x} = \{z\in Z : z< x\} \text{ and } Z_2^{x} = \{z\in Z : x\text{ and }z\text{ are incomparable}\}.
\]
Then $Z_{x}$ is the union of the disjoint sets $Z_1^{x}$ and $Z_2^{x}$. Since $X$ is image-finite, we may assume without loss of generality that $Z_{1}^{x}$ is a set of immediate predecessors of $x$.

Consider $z\in Z_2^{x}$.  Since we assumed that ${\downarrow}x\cap{\downarrow}z\neq \emptyset$, there is an element in ${\downarrow}x\cap{\downarrow}z$  of the form $\bot_{\{x,z\}}$. Let $U_{\bot_{\{x,z\}}}$ be a member of $\mathcal K$ containing $\bot_{\{x,z\}}$. By (\ref{item::top cond for bot x y}), there exists a finite subset $Q_{x,z}$ of ${\downarrow}x\cap{\downarrow}z$ such that $({\downarrow}x\cap{\downarrow}z)\smallsetminus{\downarrow}Q_{x,z}\subseteq U_{\bot_{\{x,z\}}}$.  Since $X$ is image-finite, we may assume without loss of generality that $Q_{x, z}$ is a set of immediate predecessors of $x$.

Define
\begin{equation}\label{Eq:def-of-P}
{\mathcal K}_x = \{U_x\}\cup\{U_{\bot_{\{x,z\}}} : z\in Z_2^{x}\} \text{ and } P_x = Z_1^{x}\cup\bigcup \{Q_{x,z} : z\in Z_2^{x}\}.
\end{equation}
Notice that $\mathcal{K}_{x}$ is a finite subset of $\mathcal{K}$ and $P_{x}$ is a finite subset of the immediate predecessors of $x$. Therefore, it only remains to prove that ${\downarrow}x\smallsetminus {\downarrow}P_x\subseteq \bigcup {\mathcal K}_x$, i.e., that ${\downarrow}x \smallsetminus \bigcup {\mathcal K}_x \subseteq {\downarrow}P_x$. To this end, consider $y \in {\downarrow}x \smallsetminus \bigcup {\mathcal K}_x$. Since $y\not\in U_x$ and $y \leq x$, we have $y\in{\downarrow}Z_x = {\downarrow}(Z_{1}^{x} \cup Z_{2}^{x})$. If $y \in{\downarrow}Z_1^{x}$, then $y\in{\downarrow}P_x$ as desired. So suppose that $y\leq z$ for some $z\in Z_2^{x}$. Then $y\in{\downarrow }x\cap{\downarrow}z$. Together with $y\not\in U_{\bot_{\{x,z\}}}$ this yields  $y\in{\downarrow}Q_{x,z}\subseteq {\downarrow}P_x$.
\end{proof}

We define recursively a sequence of subsets $S_n$ of $Y$ as follows:
\[
S_{1} = Y_{1} \text{ and } S_{k+1} = \bigcup\{P_y :y\in S_k \} \text{ for every integer }k \geq 1,
\]
where $Y_{1}$ is defined in (\ref{Eq:definition-of-Y}) and $P_y$ in (\ref{Eq:def-of-P}). 
Then each $S_{k}$ is a set of elements of depth $k$ in $X$. For every positive integer $k$, set
\begin{equation}\label{Eq:the-def-of-the-sets-Sk}
S_{\leq k} = S_{1} \cup \dots \cup S_{k}.
\end{equation}
By Claim \ref{claim:: building finite branching S}, the sets $S_{1}, \dots, S_k$ are finite, whence we obtain:

\begin{Claim}
\label{claim:: S_leq n is finite}
For every positive integer $k$, the set $S_{\leq k}$ is finite.
\end{Claim}

Define
\begin{align}\label{Eq:order-4}
\begin{split}
S & = \bigcup\{S_k : k\text{ is a positive integer}\},\\
S^+ & = S \cup \{\bot_A : A\in {\mathcal A}_{\infty}\text{ and }A\cap S\text{ is infinite}\}.
\end{split}
\end{align}

For every $A\in{\mathcal A}_\infty$ such that the set $A\cap S$ is infinite, choose an element $U_{\bot_A} \in \mathcal{K}$ containing $\bot_A$. By (\ref{item::top cond for bot A}), there exists $x\in A$ such that ${\downarrow}x\subseteq U_{\bot_A}$. As $\langle A, \leq \rangle$ is an image-finite diamond sequence and $x \in A$, by Proposition \ref{Prop:shapes} the set $A \smallsetminus {\downarrow}x$ is finite. Since $A \cap S$ is infinite, this implies that $A \cap S \cap {\downarrow}x \ne \emptyset$. Therefore, there is $x_{A} \in S \cap A$ such that $x_A\leq x$. Thus,
\begin{equation}\label{Eq:where-UA-is-defined}
{\downarrow}x_{A}\subseteq U_{\bot_A}.
\end{equation}
 Bearing this in mind, set
\begin{equation}\label{Eq:where-xA-is-defined}
B = \{x_A : A\in{\mathcal A}_\infty\text{ and }A\cap S\text{ is infinite}\} \text{ and } S^- = S\smallsetminus {\downarrow} B.
\end{equation}

We show that  $S^-$ is finite. For this we use the classical  Brouwer's fan theorem;  see for instance \cite[Thm.\ 3.3.20]{Da04b}. We recall that a poset is a \emph{fan} if it is rooted, its principal downsets are finite chains, and every element has only finitely many immediate successors.

\begin{Theorem}[Brouwer's Fan Theorem]
\label{theorem:: Fan theorem}
 If all chains in a fan $F$ are finite, then their length is bounded by some positive integer.
\end{Theorem}

We  turn$\langle S^-, \leq \rangle$ into a tree by means of a variant of the standard technique in modal logic known as \textit{unravelling}; see for instance \cite[Thm.\ 2.19]{ChZa97}. Consider the following set of finite sequences
\begin{align*}
T = \{ \langle x_1,x_2\ldots,x_n \rangle \;\colon &n\text{ is a nonnegative integer},\text{ }x_1 \in S_{1} \cap S^-,\ldots,x_n\in S_{n} \cap S^-,\\
& \text{ and } x_{n} < \dots < x_{2} < x_{1} \}.
\end{align*}
We assume that the empty sequence, denoted by $\langle \rangle$, belongs to $T$.\footnote{Notice that our definition differs from the standard one in that our order is $x_{n} < \dots  < x_{1}$ as opposed to the standard order $x_{1} < \dots  < x_{n}$.}

Notice  that  since the elements of each $S_{k}$ have depth exactly $k$,  the assumptions  that $x_{n} < \dots < x_{2} < x_{1}$ and $x_1 \in S_{1},\ldots,x_n\in S_{n}$  imply  that $x_{i}$ is an immediate successor of  $x_{i+1}$ for every $i < n$. We define an order relation $\sqsubseteq$ on $T$ by setting for every $s,t\in T$ that
\begin{align*}
s \sqsubseteq t \Longleftrightarrow s \text{ is an initial segment of }t.
\end{align*}

\begin{Claim}\label{Claim:F-is-a-fan}
The poset $T \coloneqq \langle T, \sqsubseteq \rangle$ is a fan.
\end{Claim}

\begin{proof}
 It is well known that  $T$ is a tree with root $\langle \rangle$ (see, e.g., \cite[Thm.\ 2.19]{ChZa97}). To see  that every element of $T$ has only finitely many immediate successors, consider a sequence $\langle x_{1}, \dots, x_{n}\rangle \in T$. The immediate successors of $\langle x_{1}, \dots, x_{n}\rangle$ have the form $\langle x_{1}, \dots, x_{n}, y\rangle$ where $y \in S_{n+1}$ and $y < x_{n}$. Since the set $S_{n+1}$ is finite, we conclude that $\langle x_{1}, \dots, x_{n}\rangle$ has only finitely many immediate successors.
\end{proof}

\begin{Claim}
\label{claim:: S - down-set of B is finite}
The set $S^-$ is finite.
\end{Claim}

\begin{proof}
Since $\vert S^- \vert\leq \vert T\vert$, it  is enough to show that $T$ is finite.
In view of Brouwer's Fan Theorem and Claim \ref{Claim:F-is-a-fan}, it suffices to prove that $T$  does not have infinite chains.

Suppose, with a view to contradiction, that there is an infinite chain $C$ in $T$. Define
\[
D = \{x\in S^- : \text{there  exists  } \langle x_1,\ldots,x_n \rangle \in C\text{ such that } x=x_n\}.
\]
Observe that $C$ is a set of finite sequences without repetitions of elements of $D$.  Since $C$ is infinite, so is $D$.  We show \color{black} that $D$ is a chain in $\langle S^-, \leq \rangle$. Let $y,z\in D$.  Then there are sequences $s=(s_1,\ldots,s_n)$ and $t=(t_1,\ldots,t_k)$ in $C$ such that $y=s_n$ and $z=t_k$. Since $C$ is a chain, either $s\sqsubseteq t$ or $t\sqsubseteq s$ in $T$. By symmetry we may assume that $s \sqsubseteq t$. Then $n\leq k$ and $s=(t_1,\ldots,t_n)$. Consequently, $y=s_n = t_n \geq \cdots \geq t_k=z$, as desired. Thus, $D$ is an infinite chain in $\langle S^{-}, \leq \rangle$.

Define
\[
A = \bigcup\{{\uparrow}y : y\in D\},
\]
where the upsets are computed in $X$. We show that $A$ is an infinite  maximal  diamond sequence in $X$.
First, $A$ is infinite  since $D\subseteq A$. Because $A$ is a subset of $X$ which is a diamond system, it is a diamond system itself. Moreover, $A$ is downward directed since $D$ is a chain. Hence, $A$ is a diamond sequence in $X$.

To prove that $A$ is maximal, consider  a  diamond sequence $A'$ in $X$ such that $A\subseteq A'$.  Let $x\in A'$. Suppose  the depth of $x$ is $n$ and let  $y$ be an element of $A$ of depth ${n+1}$  (its existence follows from the fact that $D$ is infinite). Since $x, y \in A'$ and $A'$ is downward directed, there is $z \in A'$ such that $z \leq x, y$. Thus, there is $v \in X$ of depth $n+1$ such that $z \leq v \leq x$. If $v = y$, then $y \leq x$, whence $x \in A$ (as $y \in A$ and $A$ is an upset of $X$). Otherwise, having the same depth, $v$ and $y$ are incomparable. Therefore, applying the three point rule to $z \leq v, y, x$ and $v \leq x$, we get $y \leq x$.  Thus,  $x \in A$, whence $A' \subseteq A$. We conclude that $A$ is an infinite maximal diamond sequence in $X$, i.e., $A \in \mathcal{A}_{\infty}$. In particular, this implies  that  $X^{+}$ contains an element of the form $\bot_A$.

Now, since $D \subseteq A$ is infinite and $D \subseteq S$, also  $A \cap S$ is infinite. Then the set $B$, defined in (\ref{Eq:where-xA-is-defined}), contains an element of the form $x_{A}$. By the definition of $A$, there exists $y\in D$ such that $y\leq x_A$. Therefore,  $y\in{\downarrow}B$. But this contradicts the fact that $D\subseteq S^- = S \smallsetminus {\downarrow}B$.
\end{proof}

Let $B'$ be the set of maximal elements of the set $\langle B, \leq \rangle$, defined in (\ref{Eq:where-xA-is-defined}). Since $X$ is image-finite, ${\downarrow}B'={\downarrow}B$.

\begin{Claim}
\label{claim::B' is finite}
The set $B'$ is finite.
\end{Claim}

\begin{proof}
Recall from (\ref{Eq:where-xA-is-defined}) that $B \subseteq S$ and $S^{-} = S \smallsetminus {\downarrow}B$. We show  that every element of $B'$ is either in $S_{1}$ or  has  an immediate successor in $S^{-}$.  Let  $x \in B' \smallsetminus S_{1}$. Since $B' \subseteq S$, we get $x \in S_{k}$ for some integer $k > 1$.  By definition of $S_{k}$, there is an immediate successor $y \in S_{k-1}$ of $x$. Notice that $x$ is maximal in  $B$ and $x < y$.  Consequently, $y \notin {\downarrow}B$, whence $y \in S_{k-1} \smallsetminus {\downarrow}B \subseteq S^{-}$. Thus, $x$ has an immediate successor in $S^{-}$. 

By Claim \ref{claim:: S - down-set of B is finite}, there is a positive integer $k$ such that $S^-\subseteq S_{\leq k}$. Recall that $S_{n}$ is a set of elements of  depth  $n$ in $X$, for every positive integer $n$. Therefore, the set of elements of $S$ with an immediate successor in $S^{-} \subseteq S_{\leq k}$ is  contained in  $S_{\leq k+1}$. As every element  of $B' \subseteq S$ is either in $S_{1}$ or has an immediate successor in $S^{-}$, we conclude that $B'\subseteq S_{\leq{k+1}}$. By Claim \ref{claim:: S_leq n is finite}, the set $S_{\leq{k+1}}$ is finite. Thus, so is $B'$. 
\end{proof}

\begin{Claim}\label{Fact:compactness-of-X+}
The space $X^+$ is compact.
\end{Claim}

\begin{proof}
Recalling that 
$\mathcal K$ is an open cover of $X^+$, define 
\[
{\mathcal K}' = \{U_\top\}\cup\bigcup \{{\mathcal K}_x : x\in S^-\}\cup\{U_{\bot_A} : x_A\in B'\},
\]
where $U_{\top}$  is  defined before Claim \ref{claim:: building finite branching S}, each $\mathcal{K}_{x}$ is defined in (\ref{Eq:def-of-P}), and each $U_{\bot_{A}}$ is defined in the paragraph preceding (\ref{Eq:where-UA-is-defined}). It follows from the definition that $\mathcal{K}' \subseteq \mathcal{K}$. By Claims \ref{claim:: building finite branching S}, \ref{claim:: S - down-set of B is finite} and \ref{claim::B' is finite}, the set ${\mathcal K}'$ is finite. Thus, it only remains to prove that $\mathcal{K}'$ covers $X^{+}$. Let $x\in X^+$. We have two cases: either $x \in S^+$ or $x \notin S^+$.
 
First suppose that $x \in S^+$. If $x \in S^-$, then $x\in U_x$ 
by the definition of $\mathcal{K}_{x}$ in (\ref{Eq:def-of-P}). Thus, $x \in U_{x} \in \mathcal{K}_{x} \subseteq \mathcal{K}'$. Next suppose that $x \in S^{+} \smallsetminus S^{-}$. Then \color{black} 
\begin{equation}\label{Eq:one-of-the-last-inclusions}
x\in  S^{+} \smallsetminus S^{-} \subseteq \{ \bot_{A} \colon A \in \mathcal{A}_{\infty} \text{ and }A \cap S \text{ is infinite} \} \cup {\downarrow}B \subseteq {\downarrow}B = {\downarrow}B'.
\end{equation}
The above inclusions are justified as follows. The first one follows from the definitions of $S^{+}$ and  $S^{-}$. 
The second one from the fact that for every $A \in \mathcal{A}_{\infty}$ such that $A \cap S$ is infinite, the element $x_{A} \in B$ (defined in the paragraph before (\ref{Eq:where-xA-is-defined})) belongs to $A$, whence $\bot_{A} \leq x_{A} \in B$. The last one  is  justified right before Claim \ref{claim::B' is finite}. Finally, by (\ref{Eq:one-of-the-last-inclusions}) there is $x_{A} \in B'$ such that $x \leq x_{A}$. Together with the definition of $U_{\bot_{A}}$, this implies that $x \in {\downarrow}x_{A} \subseteq U_{\bot_{A}} \in \mathcal{K}'$.

Next we consider the case where $x \notin S^+$. We have three subcases:

\benroman
\item\label{item:compactness--item-1} ${\uparrow}x\cap S=\emptyset$;
\item\label{item:compactness--item-2} ${\uparrow}x\cap S$ is infinite;
\item\label{item:compactness--item-3} ${\uparrow}x\cap S$ is nonempty and finite.
\eroman

(\ref{item:compactness--item-1}):  Since  $Y_1= S_{1} \subseteq S$, we get $x\not\in {\downarrow}Y_1$. As $Y_{1}$ is the set of maximal elements of ${\downarrow}Z$ and $X$ is image-finite, this implies that $x \notin {\downarrow}Z$. Thus, by (\ref{Eq:the-definition-of-the-set-Z-}),  $x \in U_{\top} \in \mathcal{K}'$. 

(\ref{item:compactness--item-2}):
Since $X$ is image-finite, necessarily $x \in X^{+} \smallsetminus X$. Thus, there is $A\in {\mathcal A}_\infty$ such that $x=\bot_A$. As we pointed out after the definition of $X^{+}$, we have ${\uparrow}x = \{ \bot_{A} \} \cup A$. As $\bot_{A} \notin S$ (by definition of $S$), this implies that $A\cap S={\uparrow}x\cap S$. Therefore, by assumption, $A \cap S$ is infinite, and hence $x = \bot_A\in S^+$, a contradiction with the assumption that $x \notin S^+$.

(\ref{item:compactness--item-3}): Let $y$ be a minimal element in ${\uparrow}x\cap S$. Since $x \notin S^+$ and $S \subseteq S^+$, we have $x \notin S$,  whence $x<y$. As $y \in S$, either  $y \in {\downarrow}B$ or $y \in S \smallsetminus  {\downarrow}B = S^{-}$. If $y \in {\downarrow}B$, then there is $x_{A} \in B'$ such that $y \leq x_{A}$ (since ${\downarrow}B = {\downarrow}B'$). By (\ref{Eq:where-UA-is-defined}) and $x \leq y \leq x_{A}$, we  obtain  $x \in U_{\bot_{A}} \in \mathcal{K}'$, as desired. Next suppose that  $y\in S^-$. Suppose, with a view to contradiction, that $P_{y} \cap {\uparrow}x \ne \emptyset$, where $P_{y}$ is defined in (\ref{Eq:def-of-P}). Then there is  $x \in P_{y}$ with $x \leq z$. By Claim \ref{claim:: building finite branching S}, $z$ is an immediate predecessor of $y$. Therefore, $x \leq z < y$. Moreover, from $y \in S$ it follows that there is a positive integer $k$ such that $y \in S_{k}$. By definition of $S_{k+1}$ and $z \in P_{y}$, we have $z \in S_{k+1} \subseteq S$. Together with $x \leq z < y$, this contradicts the minimality of $y$ in ${\uparrow}x\cap S$. Thus, $P_{y} \cap {\uparrow}x = \emptyset$, i.e., $x\not\in {\downarrow}P_y$. Consequently, $x \in {\downarrow}y \smallsetminus {\downarrow}P_y$. By Claim \ref{claim:: building finite branching S}, $x \in \bigcup\mathcal{K}_{y}$. Since $y \in S^{-}$, we have $\bigcup\mathcal{K}_{y} \subseteq \bigcup \mathcal{K}'$. We conclude that $x \in \bigcup \mathcal{K}'$.
\end{proof}


The final step consists in proving that  the Priestley separation axiom holds in $X^{+}$. 

\begin{law}
We call \textit{distinguished} the elements of $X^{+}$ of the form $\top$, $\bot_{\{x,y\}}$, and $\bot_A$. The remaining elements of $X^{+}$ will be  called \emph{nondistinguished}.
\end{law}

\begin{Claim}
\label{claim::principal down-sets are open}
The following conditions hold for every $u, v \in X^{+}:$
\benroman
\item\label{item::principal down-sets are open-1} The set ${\downarrow}u\smallsetminus \{u\}$ is open in $X^+$.
\item\label{item::principal down-sets are open-2} If $v$ is nondistinguished, the set ${\downarrow}v$ is  open in $X^+$.
\eroman\end{Claim}

\begin{proof}
Set $U = {\downarrow}u\smallsetminus \{u\}$ and $V = {\downarrow}v$. We verify that both $U$ and $V$ are open at once. Since $\top$ is maximal and distinguished, it does not belong to $U \cup V$. Thus, (\ref{item::top cond for top}) holds  trivially for $U$ and $V$. Moreover, since $U$ and $V$ are downsets, Condition (\ref{item::top cond for down-set x})  also holds trivially.

For (\ref{item::top cond for bot x y}), suppose that an element of the form $\bot_{\{x,y\}}$ belongs to $U$, i.e., $\bot_{\{ x, y \}} < u$. Since $X^+$ is a diamond system, it has width at most two, whence $x$ and $y$ are the only immediate successors of $\bot_{\{ x, y \}}$ in $X^{+}$. As $X$ is image-finite, ${\uparrow}\bot_{\{ x, y \}} \subseteq X$, and $\bot_{\{ x, y \}} < u$,  either $x \leq u$ or $y \leq u$. In both cases, ${\downarrow}x\cap{\downarrow}y \subseteq{\downarrow}u\smallsetminus\{u\} = U$, as desired. This arguments can be adapted to the case of $V$. For suppose that $\bot_{\{x,y\}} \in V = {\downarrow}v$. As $v$ is nondistinguished, this implies $\bot_{\{x, y \}} < v$. Thus, replicating the proof described for $U$,  we obtain that (\ref{item::top cond for bot x y}) holds for $V$ as well.

Finally, let $A \in \mathcal{A}_{\infty}$ and  $\bot_{A} \in U$. Then $\bot_{A} < u$. By definition of the order  on $X^{+}$, we have that  $u \in A$. As $A$ is infinite and ${\uparrow}u$ is finite, there  is  $x\in A$ such that $u \nleq x$. Since $A$ is downward directed (being a diamond sequence), we may assume without loss of generality that $x < u$. Then ${\downarrow}x\subseteq{\downarrow}u\smallsetminus\{u\} = U$, whence $U$ satisfies (\ref{item::top cond for bot A}). Again, this arguments can be adapted to the case of $V$ as follows. Suppose that $\bot_{A} \in V= {\downarrow}v$. Since $v$ is nondistinguished, this implies  that $\bot_{A} < v$, so we can repeat the argument for $U$.
\end{proof}

\begin{Claim}
\label{claim::deal with bot x y}
Let $x,y,u\in X$ be such that $x$ and $y$ are incomparable and ${\downarrow}x\cap{\downarrow}y\neq \emptyset$. If $\bot_{\{x,y\}}\not< u$, then either ${\downarrow}x\cap{\downarrow}y\cap{\downarrow}u=\emptyset$ or $u\in{\downarrow}x\cap{\downarrow}y$.
\end{Claim}

\begin{proof}
Recall that $x$ and $y$ are the  only immediate successors of $\bot_{\{x,y\}}$  in $X$ and hence in $X^+$, and that ${\uparrow}\bot_{\{x, y \}} \subseteq X$. 
Consequently, from $\bot_{\{x,y\}}\not< u$ it follows  that $x,y\nleq u$. Suppose there is $z\in{\downarrow}x\cap{\downarrow}y\cap{\downarrow}u$. Since $X^{+}$ has width at most two and $x, y$ are incomparable, we may assume by symmetry that $u$ is comparable with $x$. Together with $x \nleq u$, this yields $u < x$. Suppose, with a view to contradiction, that $u \nleq y$. As $y \nleq u$, this implies that $u$ and $y$ are incomparable. Since  $z \leq u, x, y$  and $u \leq x$, we can apply the three point rule to obtain $y \leq x$, a contradiction. Thus, $u \leq y$,  and so  $u \in{\downarrow}x\cap{\downarrow}y$.
\end{proof}

\begin{Claim}
\label{claim::deal with bot A}
Let $A\in \mathcal A_\infty$ and $u\in X$. If $\bot_A\not< u$,  then there is  $x\in A$ such that ${\downarrow}u\cap{\downarrow} x=\emptyset$.
\end{Claim}

\begin{proof}
Since $\bot_A\not< u$, the definition of the order relation  on $X^+$ gives that
\begin{equation}\label{Eq:the-elements-in-a-have-this-property}
a \nleq u \text{  for all }a \in A.
\end{equation}
Since ${\uparrow}u$ is finite (as $u \in X$ and $X$ is image-finite) and $A$ is infinite, there is $y\in A$ such that $u\not\leq y$. Then $u \in X$ and $y \in A$ are incomparable by (\ref{Eq:the-elements-in-a-have-this-property}). Consider $x \in A$ such that $x<y$. Its existence follows from the fact that $A$ is infinite and downward directed and ${\uparrow}y$ is finite (as $y \in A \subseteq X$ and $X$ is image-finite). Since $u$ and $y$ are incomparable and $x < y$, we get $u \nleq x$. In addition, $x \nleq u$ by (\ref{Eq:the-elements-in-a-have-this-property}). Thus, $x$ and $u$ are also incomparable. To conclude the proof, suppose that there is $z \in {\downarrow}u\cap{\downarrow} x$. As $z \leq u, x, y$ and $x, u$ are incomparable and $x\leq y$, we can apply the three point rule to obtain $u \leq y$, a contradiction.
\end{proof}

\begin{Claim}
\label{claim:: principal down-sets are clopen}
If $u\in X^+$ is nondistinguished, then ${\downarrow}u$ is clopen in $X^+$.
\end{Claim}

\begin{proof}
By Claim \ref{claim::principal down-sets are open}(\ref{item::principal down-sets are open-2}), the set ${\downarrow}u$ is open in $X^+$. Then it only remains to verify that $X^+\smallsetminus {\downarrow}u$ is also open. First, observe that Condition (\ref{item::top cond for top}) holds for $X^+\smallsetminus {\downarrow}u$ (just take $Z = \{ u \}$). To prove Condition (\ref{item::top cond for down-set x}),  let  $x\in X^+\smallsetminus {\downarrow}u$. Since $u \notin {\uparrow}x$ and ${\downarrow}x\smallsetminus{\downarrow}u\subseteq X^+\smallsetminus{\downarrow}u$, we may take $Z = \{ u \}$.

To prove Condition (\ref{item::top cond for bot x y}), suppose that $X^+\smallsetminus{\downarrow}u$ contains an element of the form $\bot_{\{x,y\}}$. Since $\bot_{\{ x, y \}} \in {\downarrow}x \cap {\downarrow}y$ and $\bot_{\{ x, y \}} \nleq u$, we  may apply Claim \ref{claim::deal with bot x y} to obtain  that either ${\downarrow}x\cap{\downarrow}y\cap{\downarrow}u=\emptyset$ or $u\in{\downarrow}x\cap{\downarrow}y$. If ${\downarrow}x\cap{\downarrow}y\cap{\downarrow}u=\emptyset$, then ${\downarrow}x\cap{\downarrow}y\subseteq X^+\smallsetminus {\downarrow}u$, while if $u\in{\downarrow}x\cap{\downarrow}y$, then we take $Z = \{ u \}$. Thus, Condition (\ref{item::top cond for bot x y}) holds in both cases.

It only remains to check (\ref{item::top cond for bot A}). Consider $A\in\mathcal A_\infty$ such that $\bot_A \in X^+\smallsetminus {\downarrow}u$. In particular, $\bot_{A} \not< u$. Since $u$ is nondistinguished by assumption and $X^{+} \smallsetminus X$ is a set of distinguished elements, $u \in X$. Therefore,  $u \in X$ and $\bot_{A} \not< u$.  Thus, we may apply Claim \ref{claim::deal with bot A} to obtain $x\in A$ such that ${\downarrow}x\subseteq X^+\smallsetminus{\downarrow}u$.
\end{proof}

\begin{Claim}
\label{claim:: 2-down-sets are clopen}
If $u$ and $v$ are incomparable elements in $X$, then the set ${\downarrow}u\cap{\downarrow}v$ is clopen in $X^+$.
\end{Claim}

\begin{proof}
Define $U = X^+\smallsetminus ({\downarrow}u\cap{\downarrow}v)$. Since $u$ and $v$ are incomparable,
\[
{\downarrow} u\cap{\downarrow} v=({\downarrow} u\smallsetminus\{u\})\cap({\downarrow} v\smallsetminus\{v\}).
\]
Therefore, in view of Claim \ref{claim::principal down-sets are open}, the set ${\downarrow}u\cap{\downarrow}v$ is open in $X^+$.

 It remains  to show that $U$ is also open.  First, observe that $X^+\smallsetminus {\downarrow}\{u,v\} \subseteq U$, whence $U$ satisfies (\ref{item::top cond for top}). To prove (\ref{item::top cond for down-set x}),  let  $x\in U$. By definition of $U$, either $x\not\leq u$ or $x\not\leq v$. By symmetry, we  may assume that  $x \nleq u$. Together with the obvious inclusion ${\downarrow}x\smallsetminus{\downarrow}u\subseteq U$, this implies that Condition (\ref{item::top cond for down-set x}) holds for $U$. To prove (\ref{item::top cond for bot x y}), suppose that $U$ contains an element of the form $\bot_{\{x,y\}}$. Then either $\bot_{\{x,y\}}\nleq u$ or $\bot_{\{x,y\}}\nleq v$. By symmetry we may assume that  $\bot_{\{x,y\}}\nleq u$, whence $\bot_{\{ x, y \}} \not < u$. Now apply Claim \ref{claim::deal with bot x y} to obtain  that either ${\downarrow}x\cap{\downarrow}y\cap{\downarrow}u=\emptyset$ or $u\in{\downarrow}x\cap{\downarrow}y$. Thus,  Condition  (\ref{item::top cond for bot x y}) follows from the inclusion ${\downarrow}x\cap{\downarrow}y\subseteq U$ in the first case and from the inclusion $({\downarrow}x\cap{\downarrow}y)\smallsetminus{\downarrow}u\subseteq U$ in the second case.

Lastly, to prove (\ref{item::top cond for bot A}), consider $A\in\mathcal A_\infty$ such that $\bot_A\in U$. By definition of $U$, either $\bot_{A} \nleq u$ or $\bot_{A}\nleq v$.  By symmetry, we may assume that $\bot_{A} \nleq u$, so  $\bot_{A}\not< u$. Moreover,  $u \in X$ by assumption.  Thus, we  may apply Claim \ref{claim::deal with bot A} to obtain  $x\in A $ such that ${\downarrow}u\cap {\downarrow}x=\emptyset$. Consequently, also ${\downarrow}u\cap{\downarrow}v \cap {\downarrow}x=\emptyset$, and hence  ${\downarrow}x \subseteq X^{+} \smallsetminus ({\downarrow}u\cap{\downarrow}v) = U$.
\end{proof}

\begin{Claim}\label{Fact:Priestley-sep-axiom}
The space $X^+$ satisfies the Priestley separation axiom.
\end{Claim}

\begin{proof}
\label{fact:: Priestley's separation}
Let $u,v\in X^+$ with $u \nleq v$. We  must  find a clopen downset $D$ of $X^+$ such that $u \notin D$ and $v\in D$. If $v$ is nondistinguished, we may take $D = {\downarrow}v$ by Claim \ref{claim:: principal down-sets are clopen}. Suppose $v$ is distinguished. We have the following cases:

\benroman
\item\label{item:Priestley-item-1} $v = \top$;
\item\label{item:Priestley-item-2} $v = \bot_{\{ x, y \}}$ for some $x, y \in X$;
\item\label{item:Priestley-item-3} $v = \bot_{A}$ for some $A \in \mathcal{A}_{\infty}$.
\eroman

(\ref{item:Priestley-item-1}): Since $X$ is image-finite, there is a maximal element $u^\circ$ in $X^+$ such that $u\leq u^\circ$. Notice that $u^\circ\neq\top$  as  $u \nleq v = \top$. Since $\top$ is maximal, this implies $v = \top \notin {\downarrow}u^{\circ}$. Being maximal in $X^{+}$ and different from $\top$, the element $u^\circ$ is nondistinguished. Thus, ${\downarrow}u^\circ$ is clopen  by Claim~\ref{claim:: principal down-sets are clopen}.  Since $u^{\circ}$ is maximal, we can apply Proposition \ref{Prop:shapes}  and obtain  that the downset ${\downarrow}u^\circ$ is also an upset of $X^{+}$. Therefore, the set $D\coloneqq X^+\smallsetminus {\downarrow}u^\circ$ is a clopen downset of $X^{+}$ such that $u \notin D$ and $v\in D$.

(\ref{item:Priestley-item-2}): If $u\notin {\downarrow}x\cap{\downarrow}y$, then we  may take $D = {\downarrow}x\cap{\downarrow}y$ by Claim \ref{claim:: 2-down-sets are clopen}. Suppose $u\in{\downarrow}x\cap{\downarrow}y$. By Proposition \ref{Prop:shapes}, for every $z\in{\downarrow}x\cap{\downarrow}y$ there is a unique immediate predecessor $z^\circ$ of $x$ and $y$ such that $z \leq z^{\circ}$.  In particular, $u^\circ$ is an immediate predecessor of $x$ and $y$. Define
\[
D = ({\downarrow} x\cap{\downarrow} y) \smallsetminus {\downarrow}u^\circ.
\]
Clearly $u\not\in D$. Moreover, since $u^{\circ}$ and $v$ are immediate predecessors of $x$ and $y$, from $v \leq u^{\circ}$ it follows that $v = u^{\circ}$. This contradicts the fact that $u \leq u^{\circ}$ and $u \nleq v$. Therefore, $v \nleq u^{\circ}$, and hence $v=\bot_{\{x,y\}}\in D$. 

 We  show that $D$ is a clopen downset of $X^{+}$. As $X^{+}$ has width at most two, its elements may have at most two immediate  successors. It follows that among immediate predecessors of $x$ and $y$ there exists just one distinguished element, namely $\bot_{\{x,y\}}$.  Since $u\not\leq v = \bot_{\{x,y\}}$, the element $u^\circ$ is nondistinguished. Thus, by Claims  \ref{claim:: principal down-sets are clopen} and \ref{claim:: 2-down-sets are clopen}, the set $D$ is clopen. Finally, to prove that $D$ is a downset, we show that
\begin{equation}\label{Eq:prove-that-D-is-downset}
z\notin D \Longleftrightarrow z^{\circ} = u^{\circ} \text{  for  every }z\in{\downarrow}x\cap{\downarrow}y.
\end{equation}
Consider $z\in{\downarrow}x\cap{\downarrow}y$. If $z^{\circ}=u^{\circ}$, then $z \in {\downarrow}u^{\circ}$ (since $z \leq z^{\circ}$), whence $z\notin D$. Conversely, if $z\notin D$, then $z \in {\downarrow}u^{\circ}$.  Therefore, $u^{\circ}$ is the unique immediate predecessor of $x$ and $y$ above $z$. Consequently, $z^{\circ} = u^{\circ}$. This establishes (\ref{Eq:prove-that-D-is-downset}). Lastly, consider $z\in D$ and $t\leq z$. Since $t \leq z \leq x, y$, we get $t^{\circ}=z^{\circ}$. Moreover,  (\ref{Eq:prove-that-D-is-downset}) and $z \in D$ yield $z^{\circ} \ne u^{\circ}$. Thus, also  $t^{\circ} \ne u^{\circ}$. With another application of (\ref{Eq:prove-that-D-is-downset}) we conclude that $t \in D$. Hence, $D$ is a downset.

(\ref{item:Priestley-item-3}): Since $u \ne v$, we have $u \ne \bot_{A}$. Because $\bot_{A}$ is the unique element whose set of strict successors is $A$, from $u \ne \bot_{A}$ it follows that there is $z \in A$ with $u\nleq z$. If $z=\top$ or $z=\bot_{\{x,y\}}$ for some $x,y\in X$, then we take the clopen downset $D$ constructed in the previous cases. Otherwise, $z$ is nondistinguished and we may take $D = {\downarrow}z$ by Claim \ref{claim:: principal down-sets are clopen}.
\end{proof}

In view of Claims \ref{Fact:downsets-of-opens-are-opens}, \ref{Fact:compactness-of-X+}, and \ref{Fact:Priestley-sep-axiom}, we conclude that $X^{+}$ is an Esakia space.
\end{proof}

\begin{Corollary}
Let $X$ be an image-finite root system. Then there  is an Esakia space $X^+$ whose underlying poset is also a root system  such that $X$ is the image-finite part of $X^+$.
\end{Corollary}

\begin{proof}
If $X$ is a root system, the poset $X^{+}$ constructed in the proof of Proposition \ref{prop:: topologies for almost forests} is also a root system. Therefore, the result follows from Proposition \ref{prop:: topologies for almost forests}.
\end{proof}

\begin{Corollary}\label{Cor:representability-for-bounded-depth}
Let $n$ be a positive integer and  $X$  a diamond system of depth $\leq n$. Then $X$ is Esakia representable.
\end{Corollary}

\begin{proof}
If $X$ is of  depth $\leq n$, the poset $X^{+}$ constructed in the proof of Proposition \ref{prop:: topologies for almost forests} coincides with $X$. Consequently, $X$ is  Esakia  representable by Proposition \ref{prop:: topologies for almost forests}.
\end{proof}

\section{Consequences}

By Theorems \ref{Thm:diamond} and \ref{Thm:main}, a variety of Heyting algebras is such that its profinite members are profinite completions if and only it it omits the finite algebras $\textup{Up}(P_1), \textup{Up}(P_2), \textup{Up}(P_3)$, and $\textup{Up}(P_4)$. Consequently, we obtain:

\begin{Theorem}\label{Thm:decidable} 
\benroman
\item[]
\item\label{item:decidability:prob:1} The problem of determining whether a finite set of equations axiomatizes a variety of Heyting algebras whose profinite members are profinite completions is decidable.
\item\label{item:decidability:prob:2} The problem of determining whether a finite set of finite Heyting algebras generates a variety whose profinite members are profinite completions is decidable.
\eroman
\end{Theorem}

\begin{proof}
(\ref{item:decidability:prob:1}): Let $\Sigma$ be a finite set of equations and $\class{V}$ the variety axiomatized by it. By Theorems \ref{Thm:diamond} and \ref{Thm:main}, the profinite members of $\class{V}$ are profinite completions if and only if $\class{V}$ validates the Jankov formulas $\mathcal{J}(P_1), \mathcal{J}(P_2), \mathcal{J}(P_3)$, and $\mathcal{J}(P_4)$. By Jankov's Lemma, this happens precisely when $\class{V}$ contains none of $\textup{Up}(P_1), \textup{Up}(P_2), \textup{Up}(P_3)$, and $\textup{Up}(P_4)$. This condition, in turn, can be decided by checking whether the equations in $\Sigma$ hold in any of the finite algebras $\textup{Up}(P_1), \textup{Up}(P_2), \textup{Up}(P_3)$, and $\textup{Up}(P_4)$.

(\ref{item:decidability:prob:2}): The variety generated by a finite set $\class{K}$ of finite Heyting algebras  is such that its profinite members are profinite completions precisely when $\class{K}$ validates $\mathcal{J}(P_1), \mathcal{J}(P_2), \mathcal{J}(P_3)$, and $\mathcal{J}(P_4)$. Since $\class{K}$ is a finite set of finite algebras, the latter property is decidable. 
\end{proof}

We next provide some properties of diamond Heyting algebras. Recall that a variety is \textit{finitely based} if it can be axiomatized by finitely many equations, and that a class of (similar) algebras is a \textit{quasivariety} if it is closed under the formation of reduced products, subalgebras, and isomorphic copies, see, e.g., \cite[Sec.\ 1.5]{Go98a}. A variety is called \textit{primitive} if all its subquasivarieties are varieties \cite[Sec.\ 5.1.4]{Go98a}. 

\begin{Theorem}\label{Thm:properties-of-DHAs}
Varieties of diamond Heyting algebras are locally finite, primitive, and finitely based. Moreover, there are only countably many of them.
\end{Theorem}

\begin{proof}
That varieties of diamond Heyting algebras are locally finite follows from Theorem~\ref{Thm:CHA-locally-finite} and the fact that $\class{DHA} \subseteq \class{CHA}$. In \cite{Citkin78a,Citkin78a-paper} it is proved that a variety of Heyting algebras is primitive if and only if it excludes the algebras of upsets of the rooted posets in Figure~\ref{fig:forbidden ordered-setsHSP}. (The original proof in \cite{Citkin78a-paper} is in Russian. For an English version of the proof, based on Esakia duality, see \cite{BezMor19}.)

\begin{figure}[h]
	\begin{tabular}{ccccccc}
		\\
		\begin{tikzpicture}
		\tikzstyle{point} = [shape=circle, thick, draw=black, fill=black , scale=0.35]
		\node[label=below:{$P_1$}]  (label) at (0,-0.1)  {};
		\node  (0) at (0,0) [point] {};
		\node (v1) at (-.833,1) [point] {};
		\node  (v2) at (0,1) [point] {};
		\node  (v3) at (.833,1) [point] {};
		\node (1) at (0,2) [point] {};
		
		\draw  (0) -- (v1) -- (1) -- (v2) -- (0) -- (v3) -- (1);
		\end{tikzpicture}
		
		\;&\;
		\begin{tikzpicture}
		\tikzstyle{point} = [shape=circle, thick, draw=black, fill=black , scale=0.35]
		\node[label=below:{$P_2$}]  (label) at (0,-0.1)  {};
		\node  (0) at (0,0) [point] {};
		\node (v1) at (-0.833,1) [point] {};
		\node  (v2) at (.7,0.6) [point] {};
		\node  (v3) at (.7,1.4) [point] {};
		\node (1) at (0,2) [point] {};
		
		\draw   (0) -- (v1) -- (1) --  (v3) -- (v2) -- (0);
		\end{tikzpicture}
		\;&\;
		\begin{tikzpicture}
		 \tikzstyle{point} = [shape=circle, thick, draw=black, fill=black , scale=0.35]
		\node[label=below:{$P_5$}]  (label) at (0,-0.1)  {};
		\node  (0) at (0,0) [point] {};
		\node (v1) at (-0.833,1) [point] {};
		\node  (v2) at (0.833,1) [point] {};
		\node (1) at (-0.833,2) [point] {};
		
		\draw  (0) -- (v1) -- (1)  (v2) -- (0) ;
	    \end{tikzpicture}
		\;&\;
		\begin{tikzpicture}
		\tikzstyle{point} = [shape=circle, thick, draw=black, fill=black , scale=0.35]
		\node[label=below:{$F_3$}]  (label) at (0,-0.1)  {};
		
		\node  (v1) at (-0.833,1.5) [point] {};
		\node  (0) at (0,0.5) [point] {};
		\node  (v3) at (0,1.5) [point] {};
		\node (v2) at (0.833,1.5) [point] {};
		
		\draw   (v1) -- (0) -- (v3) (0) -- (v2);
		\end{tikzpicture}
        \;&\;
        \begin{tikzpicture}
        \tikzstyle{point} = [shape=circle, thick, draw=black, fill=black , scale=0.35]
       
        \node[label=below:{$P_7$}]  (label) at (0,-0.1)  {}; 
        \node  (v1) at (-0.833,1) [point] {};
        \node (0) at (0,0) [point] {};
        \node  (v3) at (0.833,1) [point] {};
        \node (v2) at (0,2) [point] {};
        \node  (d3) at (1.666,2) [point] {};
        \node  (d1) at (-1.666,2) [point] {};
        
        \draw  (v2) -- (v1) -- (0) -- (v3)  -- (v2)  (d1) -- (v1) (v3) -- (d3);
        \end{tikzpicture}		
	\end{tabular}
	\caption{The posets $P_1$, $P_2$, $P_5$, $F_3$ and $P_7$.}
	\label{fig:forbidden ordered-setsHSP}
\end{figure}
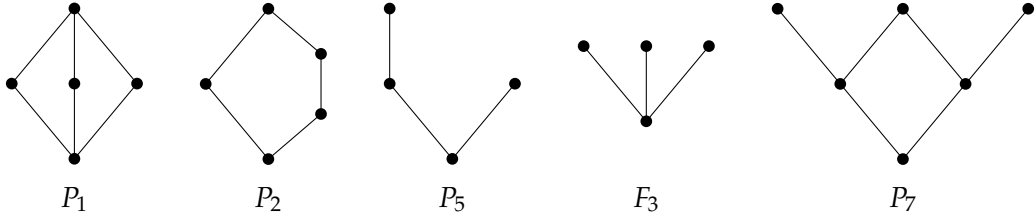

From the fact that varieties of diamond Heyting algebras omit $\textup{Up}(P_1), \textup{Up}(P_2), \textup{Up}(P_3)$, and $\textup{Up}(P_4)$,  it easily follows that they also omit the algebras of upsets of the posets $P_1, P_2, P_5, F_3$, and $P_7$. Consequently, varieties of diamond Heyting algebras are primitive. As it is known that there are only countably many primitive varieties of Heyting algebras and that all of them are finitely based (see for instance \cite[Thm.~9.1]{BezMor19}), we conclude that the same holds for varieties of diamond Heyting algebras. It only remains to prove that there are infinitely many varieties of diamond Heyting algebras. But this follows, for instance, from the fact that there are infinitely many varieties of G\"odel algebras, which are all varieties of diamond Heyting algebras. 
\end{proof}



A variety $\class{V}$ of Heyting algebras is said to be \textit{canonical} provided that $\textup{Up}(X_{\A}) \in \class{V}$ for all $\A \in \class{V}$. The terminology comes from the fact that, given a Heyting algebra $\A$, the completion $\textup{Up}(X_{\A})$ is called the \textit{canonical extension} of $\A$  \cite{GeJo94}. Furthermore, recall that the smallest variety containing a class of similar algebras $\class{V}$ is $\HHH\SSS\PPP(\class{V})$ \cite[Thm.\ II.9.5]{BuSa00}.

\begin{Theorem}
The variety $\class{DHA}$ is canonical.
\end{Theorem}

\begin{proof}
It is well known that if $\mathcal{P}$ is an elementary class of posets, then the variety of Heyting algebras generated by $\{\textup{Up}(X) \colon X \in \mathcal{P} \}$ is canonical (see, e.g., \cite[Thm.~7.1]{Go18}). Let then $\mathcal{P}$ be the class of diamond systems. Notice that $\mathcal{P}$ is elementary because the conditions in Definition~\ref{def:diamond} can be expressed by first-order sentences, so the variety 
\[
\class{V} \coloneqq \HHH\SSS\PPP\{ \textup{Up}(X) \colon X \in \mathcal{P} \}
\]
is canonical. Thus, to conclude the proof, it suffices to show that $\class{V} = \class{DHA}$.

A straightforward adaptation of the proof of the implication (\ref{item:diamond3})$\Rightarrow$(\ref{item:diamond1}) in Theorem \ref{Thm:diamond} shows that $\textup{Up}(X) \in \class{DHA}$ for all $X \in \mathcal{P}$. Consequently, $\class{V} \subseteq \class{DHA}$. To prove the other inclusion, consider $\A \in \class{DHA}$. From Theorem \ref{Thm:diamond} it follows that $X_{\A}$ is a diamond system, whence $X_{\A} \in \mathcal{P}$. Together with the fact that $\A$ embeds into $\textup{Up}(X_{\A})$ via the map $\gamma_{\A}$, this implies $\A \in \class{V}$, as desired.
\end{proof}

Intermediate logics algebraized, in the sense of \cite{BP89}, by varieties of diamond Heyting algebras have interesting metalogical properties. First, an intermediate logic is said to be \textit{hereditarily structurally complete} if all its finitary extensions \cite{AAL-AIT-f} are \textit{structurally complete} in the sense that their admissible rules are derivable (see for instance \cite{Ry97}). On the other hand, an intermediate logic $\class{L}$ has the \textit{infinite Beth definability property} if implicit definitions can be turned explicit in $\class{L}$ (we refer to \cite{BlHoo06,MRJ18sl} for the technical details).

\begin{Theorem}\label{Thm:structural-ES}
Intermediate logics algebraized by varieties of diamond algebras are hereditarily structurally complete and have the infinite Beth definability property.
\end{Theorem}

\begin{proof}
As detailed in \cite{Ry97} and \cite{BlHoo06}, an intermediate logic is hereditarily structurally complete (resp.\ has the infinite Beth definability property) precisely when it is algebraized by a variety of Heyting algebras that is primitive (resp.\ in which epimorphisms are surjective). The fact that varieties of diamond Heyting algebras are primitive follows from Theorem \ref{Thm:properties-of-DHAs}, while the fact that epimorphisms are surjective in them is a consequence of \cite[Thm.\ 9.4]{MorWan19es} (see also \cite{BMR16}).
\end{proof}

We close this paper by drawing connection to the classical Representation Problem mentioned in the introduction.

\begin{law}
We call a variety $\class{V}$  of Heyting algebras \textit{representable} if $X$ is Esakia representable for every poset $X$ such that $\textup{Up}(X) \in \class{V}$. 
\end{law}

Recall that $\class{D}_n$ is the variety of Heyting algebras of depth $\leq n$. Our aim is to prove the following: 

\begin{Theorem}\label{Thm:representable-Esakia}
A variety $\class{V}$ of Heyting algebras is representable if and only if $\class{V} \subseteq \class{D}_n \cap \class{DHA}$ for some positive integer $n$.
\end{Theorem}

To this end, we rely on the following folklore observation.

\begin{Proposition}\label{Prop:splitting-chains}
For a variety $\class{V}$ of Heyting algebras, the following are equivalent :
\benroman
\item the members of $\class{V}$ have depth $\leq n$ for some positive integer $n$;
\item $\class{V}$ omits all infinite chains;
\item $\class{V}$ omits some chain;
\item $\class{V}$ validates $\mathcal{J}(C)$ for some finite chain $C$.
\eroman
\end{Proposition}


\begin{proof}[Proof of Theorem \ref{Thm:representable-Esakia}]
Let $Z^-$ be the poset of negative integers with the standard order. First suppose that $\class{V}$ is representable. By Proposition~\ref{Prop:forbidden-posets} and Lemma~\ref{Lem:preservation},  $\class{V}$ omits the algebras $\textup{Up}(P)$ for all $P \in\{P_1, P_2, P_3, P_4\}$. In view of Theorem \ref{Thm:diamond}, we conclude that $\class{V} \subseteq \class{DHA}$. Moreover, $\textup{Up}(Z^-) \notin \class{V}$  as $Z^-$ is not Esakia representable by Proposition \ref{Prop:Esakia-tricks}(\ref{item:chains-Esakia}). Thus, by Proposition \ref{Prop:splitting-chains} the members of $\class{V}$ have depth $\leq n$ for some positive integer $n$.

Conversely, suppose that $\class{V} \subseteq \class{D}_n \cap \class{DHA}$ for some positive integer $n$. Consider a poset $X$ such that $\textup{Up}(X) \in \class{V}$. We show that $X$ is image-finite. It suffices to prove that $X$ has bounded depth and width. Since $\textup{Up}(X)$ has depth $\leq n$, the same holds for $X$ (see, e.g., \cite[Prop.\ 2.38]{ChZa97}). Therefore, by Theorem \ref{Thm:diamond}, $\class{DHA}$ has width $\leq 2$. Thus, $\textup{Up}(X)$ also has width $\leq 2$, and hence so does $X$ (see, e.g., \cite[Prop.\ 2.39]{ChZa97}). Consequently, $X$ is image-finite. As $X$ is image-finite and $\textup{Up}(X)$ is a diamond Heyting algebra, by Corollary \ref{Cor:profinite-come-from-diamond-systems} we obtain that $X$ is a diamond system of depth $\leq n$. Thus, $X$ is Esakia  representable by Corollary \ref{Cor:representability-for-bounded-depth}.
\end{proof}

\paragraph{\bfseries Acknowledgements.}
The third author was partially supported by the Beatriz Galindo grant BEAGAL\-$18$/$00040$ and by the I+D+i research project PID$2019$-$110843$GA-I$00$, both funded by the Ministry of Science and Innovation and the Ministry of Universities of Spain. The fourth author would like to express his deep gratitude to the logic group in the Institute of Computer Science of the Czech Academy of Sciences, especially to Petr Cintula, for the hospitality and support that allowed him to conduct this research.

\bibliographystyle{plain}

\end{document}